%% file: funlk.tex
\documentclass[twoside]{article}
    \usepackage[german,english]{babel}
    \usepackage[leqno]{amstex}
    \usepackage{array,amssymb,theorem}
    \usepackage{varioref,xspace}

    \usepackage[bfhd]{prodefs} %% \usepackage[bfhd]{amsproof}
    \let\invec\vcmps
\makeatother
\begin{document}
    \input{funl0-4}
        %\input{funl00} %% Macros, Initializations, Abstract.
        %   \include{funl01}  \section{Introduction}
        %   \include{funl02}  \section{Survey}
        %    \include{funl03}  \section{Signature and Language}
        %   \include{funl04}  \section{Semantics of Functional Logic}
    \input{funl5-6}
        %    \include{funl05}
        %        \section{Syntactic Matters and the Calculus of Func. Logic}
        %    \include{funl06}  \section{Formalized Theories}
    \input{funl7-bb} %       = funl07.tex + funlk.bbl
        %    \include{funl07}
        %        \section{Obtaining a Model of a Consistent Theory}
        %    \bibliography{funlk}
\end{document}

%% file: funl0-4.tex
% Logic Eprints
%Submitted 0646 Tue Mar 14, 1995 by: a8121dab@helios.edvz.univie.ac.at (josef schoenbrunner )
%logic/schoenbrunner/fun-logic-completeness
%logic/schoenbrunner/fun-logic-completeness/funl0-4.tex
%

% FUNL0-4.TEX % Copy of funl00.tex ... funl04.tex $1
%%####################################################################
%%# funl00.tex  |95-03-12|19:31|                                     #
%%####################################################################
 %    FUNL00.TEX  %  Preamble after \begin{document}
 %--------------------  Macros
     \let\include\input
     \makeatletter
 %--------------------- Switches
 %    \suppresspagenumbs
 %--------------------- effektive Part
     \hfuzz=22pt %% overfull box < 22pt no message
     \ohnetoc
     \nocite{*}
     \makeatother
 %======================================================
 %         Head (Title Abstract Classification Keywords)
 %======================================================
 \title{\vspace*{22mm} %{20mm}%{-10mm}
     \Large\bf
     Completeness Proof of Functional Logic,
     \\ \large
     A Formalism with
     %\\
     Variable-Binding Nonlogical Symbols
     \\
     }
 \author{\large Sch\"onbrunner Josef % \quad \today \quad Version 11
     \\
     \normalsize
     Institut f\"ur Logistik der Universit\"at Wien
     \\
     \normalsize
     Universit\"atsstraae 10/11,
     A-1090 Wien (Austria)
     \\
     \normalsize
     e-mail a8121dab@@helios.edvz.univie.ac.at
 } \date{}
 \maketitle

 \begin{abstract}
 %------------------------------------- %
     We know extensions of first order logic by quantifiers of the kind
     \inanf{there are uncountable many ...}, \inanf{most ...}\
     with new axioms and appropriate semantics. %%
     Related are operations such as \inanf{set of x, such that ...},
     Hilbert's $\varepsilon$-operator, Churche's $\lambda$-notation,
     minimization and similar ones, which also bind a variable
     within some expression, the meaning of which is however partly
     defined by a translation into the language of
     first order logic. In this paper a generalization is presented
     that comprises arbitrary variable-binding symbols as non-logical
     operations. The axiomatic extension is determined by new
     equality-axioms; models allocate functionals to
     variable-binding symbols.  The completeness of this system of the so
     called {\em functional logic of 1st order}\ will be
     proved.
 \end{abstract}

 \vspace{1\baselineskip} %% {1ex}

 \ifmittoc\tableofcontents\fi
 % -----------------------------------------------+
%
%%####################################################################
%%# funl01.tex  |95-03-12|09:46|                                     #
%%####################################################################
 %                    FUNL01.TEX                  %
 % -----------------------------------------------+
 % Eigentlicher Text zum Thema Funktionalen Logik %
 % -----------------------------------------------+
 \section{Introduction}
 %----------------
 \begingroup %>>>>>>>>>>>>>>>>
 \renewcommand{\dots}{...}

 Functional logic is a generalization of first order predicate logic with
 different kinds of objects by adding the following new features:

 {% -----------------------begin lokale defs ---------------------------%
 \newcommand{\vecalfa}{\alpha_1,...,\alpha_n}

 1.
     The division of expressions into the categories of sentences and
     individuals (i.e. {\em formulas } and {\em terms}) is weakened
     as with a differentiation of sorts of terms formulas shall also
     be treated as a sort. Thus the classification of the symbolic entities
     into {\em logical connectives}, {\em predicate symbols}, {\em function
     symbols } loses its significance, as the membership to one of it depends
     only on its signature (i.e. number and sorts of the argument-places and
     sort of the resulting expression). The sentential sort (formulas)
     retains its special role and will be refered to as $\prop$. Thus the
     signature of a binary connective is `$\prop(\prop,\prop)$', that of a
     $n$-ary predicate symbol `$\prop(\vecalfa)$', that of a $n$-ary function
     symbol `$\gamma(\vecalfa)$' and that of a constant symbol `$\gamma$',
     if each $\alpha_i$ and `$\gamma$' are sorts.
     Not to be found in predicate logic are symbolic entities
     whose argument-places are mixed, partly of sort $\prop$ and partly of
     another object-sort. These do not fit into any of the
     categories of {\em logical connectives}, {\em predicate symbols} or
     {\em function symbols } mentioned above. An example is the expression
     `$?(E,a,b)$' denoting an object ``$a$\ if $E$, $b$ otherwise'', which is
     built up by a symbolic entity `?' of the signature
     `$\alpha(\prop,\alpha,\alpha)$'.

 2.
     In a formalized theory of predicate logic
     expressions such as $\clabst{x}{\E}$, $\clabst{x\in M}{\E}$,
     $\iota x({\E})$, $\varepsilon x(\E)$, $\mu x(\E)$,
     $\mathop{\mu\,x}\limits_{x<b}(\E)$, $\int_a^b e \cdot dx$ are
     characterized only by an external rule of translation into the language
     of the theory. In functional logic, however, such expressions can be
     generated internally by symbolic entities that bind variables.
     This is the essential extension of this formalism.
 \medskip\par
 In standardized symbolisation a symbolic entity `op' of the resulting sort
     $\gamma$ with $k$ argument places of signature%
 \begin{math}
     (\alfa_i \mathbin,\sqv\beta_i)\komma
      \sqv\beta_i=(\beta\indij)_\berj
     {\scriptstyle(1\le i\le k)}
 \end{math}
 is linked with the generation rule by which
 \\[1ex]
 \noindent
 \halign{\hskip2em$#$&\quad if #\hfil\cr
     \oblt{op}        & $k=0$\quad (constant or variable) -- or --\cr
     \obl{\opexpr qa}& $k>0$\cr
 }
 }% ----------------- end lokale defs --------------------------------------%

 \noindent
 is an expression of sort $\gamma$, if each $a_i$ is an expression of sort
 $\alfa_i$ and each $\defsqidetj q$ is a sequence of variables of
 sorts-sequence $\defsqidetj\beta$.
 The case $r_i=0$ means that the optional part, which is written as
 $\mathop{\hbox{$[\mathinner{\ldotp\ldotp\ldotp\ldotp\ldotp}]$}}%
 \limits_{\smash{if\ r_i>0}}$, is to be dropped.
 This case applies to {\em logical connectives},
 {\em predicate symbols } and {\em function symbols } in all
 argument-places $i$, only quantifiers have a signature
 \mbox{`$\prop((\alfa):\prop)$' } with $r_1=1$.
 Putting the template
 $\mathop{\hbox{$[\mathinner{\ldotp\ldotp\ldotp\ldotp\ldotp}]$}}%
 \limits_{\smash{if\ r_i>0}
 }$
 around something is used to consider both cases $r_i=0$\
 as well as $r_i>0$. If $r_i>0$ the brackets can be erased and
 `$\mathop{[(\sqv q_i):]}\limits_{\smash{if\ r_i>0}} a_i$'
 stands for `$(\sqv q_i):a_i$', which is an abbreviation of
 $\obl{(q_{i\,1},\ldots,q_{i\,r_i}):a_i}$. $\obl{q_{i,j}}$
 are the binding variables to $\obl{a_i}$.

 \noindent
 {%----------- lokale \def's:
     \def\apl{\text{Apl}} \def\PR{\rmp{PR}} \def\felse{?}
     \def\pred{\text{pre}} \def\nmb#1: {\quad{#1}.\enspace}

 {\bf Examples:} \qquad 1.\quad
 The extension of a formal Peano-system by axioms like
 $$
     \apl_k((x_1,\ldots,x_k):e,a_1,\ldots,a_k)=%
     \Subs{e}{x_1\zweildt x_k}{a_1\zweildt a_k} ,\hfill
 $$
 $$
     \PR(a,(y,z):b,n)=%
     \felse(n\mathord=0,\>a,\>\apl_2((y,z):b,\>n-1,%
     \>\PR(a,(y,z):b,n-1)))
 $$
 \noindent %\item[]\stepcounter{enumi}  %\nmb1:
 allows the representation of each primitive recursive function by
 a single term.
 (In the above schemes of axioms $e,a,b,a_i$ range over arbitrary terms
 and $n$ is a number variable.)
 If all free variables of $a$ and $b$ which are not members of
 $\encurs{\obl y, \obl z}$ are in $\encurs{\obl{u_1},...,\obl{u_m}}$,
 then the term $\obl{\PR(a,(y,z):b,n)}$ can be associated with a
 $m+1-$ary function of arguments $u_1,...,u_m,n$, defined by
 primitive recursion from {\em base-function }
 $\enangle{u_1,...,u_m} \mapsto a$ and
 {\em iteration-function } $\enangle{u_1,...,u_m,y,z} \mapsto b$.

 2.
     Quantifiers to variables of different sorts must be distinguished,
     the signature of $\obl{\GQ^{\alpha}}$ is $\obl{\prop((\alfa):\prop)}$.
     In standardized manner, a formula $\obl{\AQ{x^{\alfa}} E}$ would be\\
     $\obl{\GQ^{\alpha}((x^{\alfa}):E)}$.

 3. %\item %\nmb3:
     A standardized version of expressing
     ``the least $x$ less than $b$ such that $E$ if one exists,
       or $b$ if none exists''
     (usually symbolized by $\obl{\mathop{\mu\,x}\limits_{x<\rmp b}\,\E}$)
     is $\obl{\mathop\mu_<(\rmp b,(x):\,\E)}$, the signature of
     $\obl{\mu_<}$ being $\obl{\nu(\nu,(\nu):\prop)}$
     if $\nu$ is the sort of natural Numbers.
 } %---------- End Ber. Lokaler \def's
 \vspace{1ex}
 \par
 A standardized symbolic language and an ideal language for application
 with the same expressional ability are different.
 The first should be simple in order to avoid unnecessary expense
 in metatheoretic treatment. With regard to application this simplicity can
 be disadvantageous. For instance in predicate calculus one symbol cannot
 be used with different signatures depending on the sorts of arguments
 it appears with. Such multiple use of a symbol became popular
 in programming languages, when looking at
 overloaded versions of procedure-names.
 Application of formal logics could profit from such a technique, too.
 For instance, consider sorts $\alpha,\beta$ and a class of
 models such that the range of $\beta$ is a substructure of the range of
 $\alpha$, if the signature of a symbol w.r.t a certain argument-place
 is of sort $\alpha$, then any term of sort $\beta$ also fits into that place.
 ``overloading of symbols'' may yield simpler axiom-schemes.
 Yet it requires change from the notion of {\em symbol } to that of
 {\em symbolic entity } (= symbol + signature).
 As a basis of meta-linguistic reference we shall take the standardized form.
 Results derived on this basis can easily be transferred into more flexible
 symbolism for practical use. Non standardized usages of writing
 such as that w.r.t. quantifiers
 shall be retained like alias clauses in our object language.
 Instead of overloading the various `$\GQ^\alpha$' into one `$\GQ$'
 and various `$\overset{\alpha}=$' to `=' we stipulate:
 `$\GQ x$' stands for `$\GQ^\alpha x$' if $\obl x \in \VAR_{\alpha}$
 and `$a=b$' stands for `$a \overset{\gamma}= b$'
 if $\obl a, \obl b \in \F_\gamma$.

 As to the logical axioms, the usual schemes of predicate calculus may be
 adapted, but {\em binding of variables \/} (significant to the axioms)
 is performed by symbols other than quantifiers, too.
 One part of the equality axioms become
 \begin{displaymath}%^{\renewcommand{\dots}{...}}
     \hspace*{-7mm}
     \ifBonnerQ ::::::::::%
     {\textstyle\GQ_{z\biind i1}...\GQ_{z\biind i{r_i}}}
     \else %::::::::::::::%
     (\forall z\biind i1)...(\forall z\biind i{r_i})\;
     \fi %::::::::::::::::%
         (a_i\vSubs{x_i}{z_i} \stackrel{\alpha_i}= b_i\vSubs{y_i}{z_i})
     \limp \opex xa \stackrel{\gamma}= \opex yb
 \end{displaymath}

 where $\Vec{x_i} \equiv \xtup{x}i{r_i}$, similarly $\Vec{y_i}, \Vec{z_i}$,
 and where $a_i\vSubs{x_i}{z_i}$ designates the expression obtained from $a_i$
 by replacing each {\em free occurrence } of $x\biind ij$
 by $z\biind ij$ (for {\small$1 \le j \le r_i$}) and
 $\opex yb$ differs from $\opex xa$ only by the $i$-th argument.
 \noindent
 Note:
 (1) if $r_i=0$, then the above sequence of universal quantifiers
     becomes empty;
 (2) If $\gamma=\prop$ is the sentential sort,
 then $\stackrel{\prop}=$ is to be identified with
 $\leftrightarrow$ (=logical equivalence).

 The main problem is introducing appropriate semantics
 to which the calculus is complete.
 Let ``op'' be a {\em symbol } with $k>0$ argument places,
 at least one of them provides binding variables
 i.e. $r_i>0$ for some $\indto ik$.
 At first consideration we suppose an interpretation-structure
 to assign to `op' the functional
 \vspace{-3pt}
 $$\textstyle
     \M(\oblt{op}):\Prod_{i:1\zweildt k}V_i\>\to\>\M_{\gamma}\>,\mkern 22mu
     V_i =
     \hbox{\footnotesize$%
     \begin{cases}
         \M_{\alpha_i}    & \text{if \enspace $r(i)=0$}
         \\
         \FigMap        & \text{if \enspace $r(i)>0$}
     \end{cases}$}
 $$
 \noindent
 ($\M_{\gamma}$ is the range of $\gamma$ and
  $\operatorname{Map}(\rmp X,\rmp Y)=\clabst{f}{f:\rmp X\to\rmp Y}=Y^X
  $).
 But this turns out to fix too much,
 as assignment only to a part of the functions of
 $\FigMap$
 will be relevant for evaluation of expressions.
 Nothing beyond that partial assignment you may expect to come out from the
 {\em syntactic information of a consistent theory}.
 To overcome this problem
 a certain restriction of the argument ranges $V_i$ will help.
 The notion of a structure $\M$ must therefore be extended
 by a new component which assigns a selected set
 $\M_{\gamma}^{\vec{\sigma}} \subseteq
  \text{Map}(\Prodim\M_{\sigma_i},\M_{\gamma})$
 to each sequence of sorts $\gamma,\vec{\sigma}$.
 The selected sets are characterized by some
 {\em closure qualities } similar to those that apply to
 the set of {(primitive-) recursive functions\/},
 for instance constant functions and projections are to be included.
 In a trivial way, however, we find an extension $\BM$ of $\M$ so that
 $\BM_{\gamma}^{\vec{\sigma}}=\text{Map}(\Prodim\M_{\sigma_i},\M_{\gamma})$
 and the {\em interpretations of expressions } by $\M$ and by $\BM$ coincide
 as well as the {\em semantic consequences } $\M\satq$ and\enspace $\BM\satq$.
 \nocite{goedel}
 To construct a model of a {\em consistent formal theory } the method
 of extension to a {\em complete Henkin Theory } as in Henkin's
 Proof of the Completeness Theorem s. \cite{henkin,shoenf} is still applicable.
 \endgroup %%%<<<<<<<<<<<<<<<<
%
%%####################################################################
%%# funl02.tex  |95-03-11|22:30|                                     #
%%####################################################################
 %         FUNL02.TEX          %
 % ------------------------------------+
     \section{Survey}
 %%------------------------\noindent--------------------------------------
 As basic structure of a $1^{\text{st}}$ order functional logic language
 we define the $\lfi-$signature. Then a standardized language is specified
 that determines the notion of an expression $\obl{e}$ of sort $\gamma$.
 This is defined inductively by a characteristic syntactic relation
 of $\obl{e}$ to a symbol `op' (the root of $\obl{e}$), argument expressions
 $\obl{a_i}$ and possibly variables $\obl{v\biind ij}$ binding $\obl{a_i}$.
 As this relation shall frequently appear as a background premise within
 definitions and proofs constantly using the same arguments
 $\obl e$,`op',`$a_i$' and $\obl{v\biind ij}$,
  we introduce the abbreviation \Syntass.
 The definition of a $\lfi-$structure is based
 on the notion of a $\lfi-$signature
 according to features discussed in the introduction.
     We shall only consider logic with {\em fixed equality } base on
     {\em normal structure semantics}.
 To derive semantics for the language from the notion of structure
 based on a signature, that is to establish an interpretation
 of the expressions (of various sorts), the usual definition as a map
 from {\em variables-assignments } to the domain of the sort the expression
 belongs to is not suitable. Instead of it now an expression $\obl{e}$
 will be evaluated according to a $\lfi-$structure $\M$
 by assigning a mapping on the set of the so called perspectives of $\obl e$
 consisting of all finite sequences of variables, such that all free variables
 of $\obl{e}$ appear within that sequence.
 Let $\gamma$ be the resulting sort of $\obl e$.
 The evaluation of $\obl e$ based on $\M$
 maps the empty sequence $\enangle{}$
 into a member of the range $\M_{\gamma}$  of $\gamma$,
 provided that $\enangle{}$ is a {\em perspective } of $\obl e$
 (i.e. if $\obl e$ has no free variable)
 and it maps a non-empty perspective
 $\enangle{\obl{u_1},\ldots,\obl{u_m}}$ of $\obl{e}$
 into a function of
 $\M_{\sigma_1}\times\ldots\times \M_{\sigma_m}\longrightarrow \M_{\gamma}$,
 if $\sigma_j$ is the sort of the variable
 $\obl{u_j} \enspace \scriptstyle(j=1,\ldots m)$.
 The definition will be {\em inductive }
 based on the background-assumption of \Syntass.

 As to the {\em axiomatization\/},
 the {\em logical axioms } differ in shape from {predicate logic }
 only a little with regard to {\em equality logic}.
 But we must also take into account an extension of some notions
 which are basic to formulate axioms of logic, namely the notions of
 free and bound variables, substitution and substitutability.
 The {\em axioms system } together with the {\em rules }
 {\sf Modus Ponens } and {\sf Generalization } establishes
 the {\em calculus of Functional Logic}. The extension of
 this calculus by individual {\em nonlogical axioms }
 is called a {\em functional logic theory}.
 A $\lfq-$structure-model of a
 {\em consistent functional logic theory }
 can be constructed
 as in predicate logic
 from an extension of that theory which inherits consistency,
 admits examples and is complete.
 ({\em admitting examples } is related to the existence of terms
     $t$ for each formula $\varphi$ with at most one free variable $x$,
     so that $\EQ{x}\,\varphi \limp \varphi\psub xt$ is a theorem;
     we associate this theorem to designate $t$ as an example,
     if $\EQ{x\,\varphi}$ is true.
 )
 In Henkin's proof this is achieved in two steps:
 The 1st extension produces a {\em theory } that admits examples
 by addition of constant symbols and
 {\em special axioms } (s. \cite{henkin,shoenf,barwise}).
 Consistency continues as this extension is {\em conservative }
 (each theorem of the extended theory,
 if restricted to the original language,
 is also provable within the original theory).
 The 2nd extension by Lindenbaum's theorem enlarges the set of
 nonlogical axioms without changing the language.
 Both extensions can easily be adapted to functional logic.
 The definition of a {\em\inanf{term structure}},
 which shall prove to be a model of the constructed extension to a
 {\em closed Henkin Theory} and hence also a model of the original theory,
 also relies on a so called {\em norm function} that assigns
 a representative to each {\em closed expression } within
 a {\em congruence class}.
 This class will be defined by the congruence relation,
 that applies to $\obl{a}$ and $\obl{b}$ iff $\obl{a=b}$ is a theorem of
 the {\em extended theory}.
 As we suppose {\em completeness } of the {\em extended theory},
 there are exactly two {\em congruence classes} of
 expressions of sort $\prop$ ; hence we choose
 the constants $\falsum$ and $\verum$
 (representing {\em true } or {\em false } respectively)
 as values of the {\em norm function} of formulae.
 Upon the set of {\em norms } (i.e values of the {\em norm function}),
 which is a subset of {\em closed expressions } to each {\em sort}
 as {\em base-range}, we then define our so called {\em term-structure}.
 The model quality of this structure will be obtained as
 an immediate consequence of a theorem (by specialization).
 The claim of this theorem is that the {\em evaluation } of an
 expression $\obl e$ by the {\em term-structure } $\CM$ is a function which
 assigns to each {\em perspective } a mapping from a cartesian product of
 certain ranges $\zweildt\CM_{\sigma_i}\zweildt$ to $\CM_{\gamma}$,
 which can be described exclusively by application of
 {\em multiple substitution } (variables by terms) from $\obl e$
 and application of the {\em norm function}.
 The validity of a formula (=\thinspace expression of sort $\prop$)
 within a model means that its interpretation maps one
 (and implicitly all) non-empty {\em perspectives } into a constant
 function of value $\M(\verum)$.
 In case of a closed formula this implies
 that the {\em empty perspective } is assigned the value $\M(\verum)$.
 If $\CM$\ takes the place of $\M$,
 $\M(\verum)$ changes into $\verum$\ ($=\CM(\verum)$).
 By applying the preceding theorem to an $\obl{e}$ of sort $\prop$
 and taking into account that {\em equality of the sort } $\prop$ and
 {\em logical equivalence} become one and the same
 ($\obl{\eqs{\prop}}=\obl{\leftrightarrow}$),
 you easily conclude the equivalence of $\obl{e}$
 {\em being valid in the term-model}
 and
 {\em being deducible in the extended theory}.
 As we refer to an {\em extension}, the restriction of
 $\CM$ to the language of the original theory is also a model
 of this theory. This confirms the {\em satisfiability } of that theory
 on the assumption of its {\em consistency}.
%
%%####################################################################
%%# funl03.tex  |95-03-12|10:40|                                     #
%%####################################################################
 %         FUNL03.TEX          %% adaptiert fuer LATEX am 17.02.88
 %---------------------------+
 %\rnotiz{\sc Funl03}
 %---------------------------+
 \begingroup %---------------------------------%
 \newcommand{\eln}{\\&}
 \newcommand{\elnr}{\\[\medskipamount]} \let\elnb\elnr
 \newcommand{\elns}{\elnr&}
 \newdimen\einzug \einzug=0em
 \def\postulat#1\par#2\par\par
 {\begin{center}\usebox{\postulbox}\end{center}
  % $$\rm#1$$ \hfil % \hang
   $$\displaylines{\hspace*{-1.4em}\rm#1\hfill
   \\[1.5\baselineskip]#2}$$
 }
 %%\def\itx#1{\hbox{\it#1}}
 %----------------------------------------------%
 \def\rmm#1{$\rm#1$} \def\sgrp#1#2{{#1#2}}
 \let\leqv\leftrightarrow \let\lleqv\longleftrightarrow
 \def\enbrace#1{\{#1\relax\}} \def\moqt#1{\enqt{\mathord{#1}}}
 \def\Junklst%
 {    \moqt\verum\kom \moqt\falsum\kom \moqt\neg\kom \moqt\rightarrow\kom
     \moqt\wedge\kom \moqt\vee\kom \moqt\leftrightarrow
 }
 %---------------------------------------------------------------------%
 \section{Signature and Language}
 %---------------------------------------------------------------------%
 \begin{defin} \noem \label{Df.Fnl}
 %%    \DKLBOX{\fnlsigs}{Df.Fnl}
     $\fnlsigs$\enspace:\enspace
     The notion of
     \inanf{$S$ is a signature of $1^{\text{st}}$-order functional logic}
     is determined by the following key-components:
 %    \\
     \emma{\SRT_S}: sorts;
     \emma{\SOP_S}: symbolic operations;
     \emma{\VSRT_S}: %%%%%%%% Teilmenge von \emma{\SRT_S} --
         sorts for which variables and quantification are provided.
     \emma{\VAR_S}: variables;
     \emma{\sig_S}: signature map,\enspace
     \emma{\sigsoop S= \stdsigx
     }\
     for \emma{\oblt{op} \in \SOP_S \cup \VAR_S}\
     characterizes \oblt{op} as a symbolic operation
     to generate expressions of sort $\gamma$ from $n$
     argument-expressions of sort $\alpha_i$, that might be bound
     by \emma{r_i}\ variables of sorts \emma{\beta\indij}.
     Significant for the notion to be defined
     is also a distinguished sort $\prop$ (the type of formulae)
     and distinguished elements of \emma{\SOP_S}\ :
     $%\begin{math}
         \moqt\verum\kom \moqt\falsum\kom \moqt\neg\kom
         \moqt\rightarrow\kom \moqt\wedge\kom
         \moqt\vee\kom \moqt\leftrightarrow\kom
         \moqt{\GQ^\alfa} \kom \moqt{\PQ^\alfa}
     $ %\end{math}
     (for each $\alfa$ of \emma{\VSRT_S}) and
     $\moqt{\eqs\alfa}$\ (for \emma{\alfa \in \SRT_S})
     with fixed values relative to \emma{\sig_S}.
     In formalized manner now we stipulate all
     characterizations of this definition as follows:
 \end{defin}

 \noindent
 $\fnlsigs \;\boldsymbol{\lleqv}\;
     \text{Conjunction of the following attributes}:
 $
 \\[-1\baselineskip]
 \begin{align*}%%^{\allowdisplaybreaks}
     \quad&%
         S = (\SRT_S,\SOP_S,\VSRT_S,\VAR_S,\sig_S)
     \skp
         \VSRT_S \incl \SRT_S \skp \VAR_S\cap \SOP_S=\emptyset
     \eln
         \sig_S\colon(\SOP_S \cup \VAR_S) \to
         \SRT_S \times \textstyle
         \bigcup\limits_m(\,\SRT^m \times {\VSRT^*}^m\,)
     \eln
         (\GQ \enqt v \in \VAR_S)\
         \sig_S\enqt v \in \VSRT_S\times\enbrace\enangle\empty^2
     \skp
         \SOP_S \cup \VAR_S \quad \text{can be well-ordered}
         \footnotemark[1]
     \eln
         (\GQ \alfa \in \VSRT_S)\
         \sparenth{big}{%
         \clabst{\enqt v \in \VAR_S}{%
             {\sig_S\enqt v=(\alfa,\emptyseq,\emptyseq)}}
         \quad \mbox{is enumerable}}
     \elns
         \prop \in \SRT_S \skp \Junklst \in \SOP_S
     \eln
     (\GQ \alfa \in \VSRT_S) \quad
         \moqt{\GQ\nolimits^\alfa}, \moqt{\PQ\nolimits^\alfa} \in \SOP_S
         \qquad\qquad
     (\GQ \alfa \in \SRT_S) \quad \moqt{\eqs\alfa} \in \SOP_S
     \elns
     \parbox{.94\textwidth}{%+++++++++++++++++++++++++++++++++++
         $\sig_S$ for the distinguished members of $\SOP_S$
         is specified by a circumscription $\ustypeS$,
         \enspace (s. auxiliary notations below)
     \enspace:} %+++++++++++++++++++++++++++++++++++++++++++++++
 \end{align*}

 \vspace{-1\baselineskip}
 \begin{displaymath}
 \begin{array}{*{5}{c|}c}
     \text{op}
     &\verum, \falsum &\lnot              &\limp,\land,\lor,\leqv
     &\GQ^\alfa, \PQ^\alfa                &\eqs\alfa
     \\ \hline
     \ustypeS\text{`op'}
     &\obl{\prop}    &\obl{(\prop)\prop}  &\obl{(\prop,\prop)\prop}
     &\obl{((\alfa)\prop)\prop}           &\obl{(\alfa,\alfa)\prop}
     \\
     &&&&\mbox{\footnotesize$\for \alpha\in\VSRT_S$}
     &\mbox{\footnotesize$\for \alpha\in\ISRT_S$}
     \\ \hline
 \end{array}
 \end{displaymath}

 \footnotetext[1]{%
     This is needed only to prove completeness independent of
     the {\em axiom of choice}.
     }

 \noindent
 \underbar{Auxilary Notations} \enspace
 (dependent components) to a given $\fnlsigs$:
 \\[\medskipamount]
 \begin{math}^{\renewcommand{\arraystretch}{1.2}
         \newcommand{\vsigalpha}{(\alfa,\emptyseq,\emptyseq)}
     }
     \begin{Array}{lrclc>{$\quad}r<{$}}
         \AQu{\alpha\in\SRT_S} &\CSOP_S^{\alpha} &=
         &\clabst{\oblt{c}\in\SOP_S}{\sig_S\oblt{c}=\vsigalpha}
         &\aliasgl \CSOP_{\alpha}
         &constants
     \\
     \AQu{\alfa \in \VSRT_S}
         &\VAR^\alfa_S &=
         &\clabst{\obl{v} \in \VAR_S}{{\sig_S\obl{v}=\vsigalpha}}
         &\aliasgl \VAR_{\alpha}
         &variables% v. Typ $\alpha$
     \end{Array}
     \\
     \begin{array}[t]{rl}
     \AQu{\vec{\sigma} = \enangle{\sigma_i}_{\indto il} \in \VSRT\stern}
         &\VAR^{\vec{\sigma}}_S \aliasgl \VAR_{\vec{\sigma}}
         = \Prod_{\indto il} \VAR_{\sigma_i} =
         \\
         &= \clabst{\vec{u}}{%--------------
             \vec{u} = \enangle{u_i}_{\indto il}
             \land
             \AQu{\indto il} \obl{u_i} \in \VAR_{\sigma_i}
             }
     \end{array}
 \end{math}
 \\[1ex]
     For $\sig_S$ we use a circumscription that is more convenient
     for application:
     \begin{align*}
         \ustypeS \colon \SOP_S &\cup \VAR_S \to
             (\SRT_S\cup \enbrace{\enqt(,\enqt,,\enqt)})\stern
         \qquad %%\\[\medskipamount]
         \AQu{\enqt{\text{op}} \in \SOP_S \cup \VAR_S}
     \\
         \sigsoop S &= \stdsigx
         \quad \text{ iff }
         \\
         \ustypeS\;\oblt{op} &=
         \begin{cases}
             \enqt\gamma                                &\falls m\mathbin=0
             \\
             \enqt{(\theta_1,\zweildt,\theta_m)\gamma}  &\falls m\mathbin>0
         \end{cases}
         %%% \qquad \text{where} \\ &\hspace{15mm} \AQu{i:1\zweildt m}
         %%% \quad \enqt{\theta_i}=
         \qquad \enqt{\theta_i}=
         \begin{cases}
             \enqt{\alfa_i}             &\falls r_i=0
             \\
             \enqt{(\beta_{i\;1},\zweildt,\beta_{i\;r_i})\alfa_i}
             &\falls r_i>0
         \end{cases}
     \end{align*}

 \begingroup %>>>>>>>>>>>>>>>>
     \renewcommand{\emma}[1]{$#1$}
 \begin{mkonv}\em
     \vspace{1.4ex}
 %%    \begin{enumerate}^{\renewcommand{\labelenumi}{\arabic{enumi}) }}
     \begin{parnum}
         %----------------------------------------------------------
     \item
          Subscript $S$ will be omitted ($\SRT$ for
         $\SRT_S, ... ,\mathop{\VAR_\alfa}
         \\
         \text{ for } \mathop{\VAR^\alfa_S}$)
         if only one $\fnlsig$ is considered.
     \item
         $\alpha,\beta,\gamma,...,\alpha_i,\beta_{ij},...$
         denote members of $\SRT$.
     \item
         $u,...,z,\enspace u_i,...,v_{ij},...$
         denote members of $\VAR$.
     \item
         Symbols with an arrow-accent refer to a finite sequence
         and  writing such symbols one after the other
         denotes the concatenation of the sequences
         (if $\vec{p}=\enangle{\tup pk}$ and $\vec{q}=\enangle{\tup ql}$
         then $\vec{p}\vec{q}=\enangle{\tup pk, \tup ql}$).
     \item
         If such a symbol e.g. $\vec{u}$
         appears inside a quoted string,
         as for instance \mbox{$\obl{\text{op}((\vec{u}):a)}$},
         it denotes the string $\obl{u_1,\ldots,u_l}$,
         which is the concatenation of each $\obl{u_i}$
         with $\obl{,}$ interspearsed.
     \item
         We shall always assume
         $\vec{\alpha}=\indtupel{\alpha}im \komma
             \VVec{\beta}=\indtupel{\vec{\beta}}im \komma
             \vec{\beta_i}=\enangle{\beta\biind ij}_{\indto j{r_i}}
         $.
     \end{parnum}
 \end{mkonv}
 \endgroup %%%<<<<<<<<<<<<<<<<
 \endgroup %-----------------------------------------------------------

 \begin{defin}\noem
     $\enspace{\parenth{big}{\F_S^\gamma}_{\gamma \in \SRT_S}}$
     \label{Sprache}
     Let $S$ be a \emma{\fnlsig}.
     The {\em standardized language } of $S$ is introduced as
     a mapping on $\SRT_S$ by stipulating for each \emma{\gamma \in \SRT_S}
     the set $\F_S^\gamma$
     of expressions of sort $\gamma$ inductively as follows
     (by above clause (1)\quad $\F_\gamma \aliasgl \F_S^\gamma$):
 \end{defin}
 \vspace{-1\baselineskip}
 \begin{displaymath}
     \enqt{e} \in \F_\gamma \leqv
     (\EQ{\oblt{op},m,\gamma,\vec{\alpha},\VVec{\beta},\vec{a},\VVec{v}})
     ( \Syntass)
 \end{displaymath}
 \noindent
 where \Syntass\ abbreviates the conjunction of the following formulae:
 \iftrue %:::::::::%
 %\\[1ex]%
 \begin{displaymath}
     \begin{array}{l}
     \obl{\op} \in \SOP \cup \VAR \spa \sigoop = \sigopx
      \\
      \vec{\alpha} \in \SRT^m
      \spa
      \VVec{\beta} \in \tprod{\indto im}{\VSRT^{r_i}}% {\VSRT^*}^m
      \spa \mbox{($\vec{\alpha},\VVec{\beta}$ rely on above I.(6))}
      \\
     \vec{a} = \enangle{\obl{a_i}}_{\indto im}
         \in \textstyle\Prod\limits_{\indto im} \F_{\alpha_i}
     \spa
     \mbox{(each $\obl{a_i} \in \F_{\alpha_i}$)}
     \qquad
     \VVec{v} = \enangle{\vec{v}_i}_{\indto im}
      \in \Prod\limits_{\indto im} \VAR_{\vec{\beta_i}}
     \\
     \AQ{i}\;
     \parenth{Big}{\vec{v_i} =
         \enangle{\obl{v\biind ij}}_{\indto j{r_i}}
     %    \land
     %    \AQ{\indto{j}{r_i}}\;\obl{v\biind ij} \in \VAR_{\beta\biind ij}
         \land
         \AQu{1\le j < k \le r_i}\obl{v\biind ij}\neq\obl{v\biind ik}
     }
     \iftrue %:::::::::::%
     \\
     \obl e =
     \begin{cases}
         \obl{\op} & \falls  m=0
         \\
         \obl{\opex va} & \falls m>0
     \end{cases}
     \else %::::::::::::::%
     \\[3pt]
     if $m=0$    \enspace  then\enspace $\obl e = \obl{\op}$ \qquad
     if $m\neq 0$\enspace  then\enspace $\obl e = \obl{\opex va}$
     \fi %::::::::::::::::%
     \end{array}
 \end{displaymath}
 \else %::::::::::::::%
 \vspace{-1ex}
 \begin{equation*}
     \begin{aligned}
         \obl{\op} &\in \SOP \cup \VAR
         \\
         \vec{\alpha} &=
         \enangle{\alpha_i}_{\indto im}
         \in \SRT^*
         \\
         \vec{a} &= \enangle{\obl{a_i}}_{\indto im}
         \in \textstyle\Prod\limits_{\indto im} \F_{\alpha_i}
     \end{aligned}
     \begin{aligned}
         \sigoop &= \sigopx
         \\
         \VVec{\beta} &=
         \enangle{\vec{\beta}_i}_{\indto im}
         \in \VSRT^{**}
         \\     \bigl(
                 &\leqv(\AQ{\indto im})\enspace \obl{a_i} \in \F_{\alpha_i}
             \bigr)
         \vphantom{\Prod\limits_{\indto im} \F_{\alpha_i}}
     \end{aligned}
 \end{equation*}
 \vspace*{-2\baselineskip}
     \ifkurz %:::::::::::%
 \begin{align*}
     {}&\VVec{v} = \enangle{\vec{v}_i}_{\indto im} =
     \enangle{\enangle{\obl{v\biind ij}}_{\indto j{r_i}}}_{\indto im}
         \komma \text{each } \obl{v\biind ij} \in \VAR_{\beta\biind ij}
         \komma \text{if } j\neq k \text{ then }
         \obl{v\biind ij} \neq \obl{v\biind ik}
         \\[3pt]
         \empty&
         \text{if }  m=0  \text{ then } \obl e = \obl{\op} \qquad
         \text{if }  m\neq 0  \text{ then } \obl e = \obl{\opex va}
 \end{align*}
     \else %::::::::::::::%
 \begin{align*}
     {}&\enangle{\VVec{v}} = \enangle{\vec{v}_i}_{\indto im} =
     \enangle{\enangle{\obl{v\biind ij}}_{\indto j{r_i}}}_{\indto im}
     \\
     \empty&(\AQ{\indto im})
     \begin{Array}^{\renewcommand{\arraystretch}{1.2}}[t]{ll}
         (\AQ{1 \le j \le r_i})
         & \obl{v\biind ij} \in \VAR_{\beta\biind ij}
         \\
         (\AQ{1 \le j < k \le r_i}) \enspace
         & \obl{v\biind ij} \neq \obl{v\biind ik}
     \end{Array}
     \\
     \empty&\obl e =
     \begin{cases}
         \obl{\op} & \falls  m=0
         \\
         \obl{\opex va} & \falls m>0
     \end{cases}
 \end{align*}
     \fi %::::::::::::::::%
 \fi %::::::::::::::::%

 This definition characterizes the expression $\obl e$ as a chain of
 symbols which is produced by a symbolic operation $\oblt{op}$ either
 exclusively (constant or variable) or together with argument-expressions
 $\obl{a_i}\quad (\indto im)$ possibly accompanied by binding variables
 $\vec{v_i}$ ($\obl{e} \in \F_{\gamma}$ is composed of smaller
 expressions $\obl{a_i} \in \F_{\alpha_i}$).

 In predicate logic binding variables $\vec{v_i}$ are only provided for
 the two quantifiers, but expressions which are built up by another symbol
 $\obl{\op}$ are either of shape $\obl{\op}$ or $\obl{\op(a_1,\zweildt,a_m)}$.
 Even in application of functional logic
 binding variables will be rare and never appear in front of more than
 one argument of a symbolic operation.
 The above definition is a prerequisite to almost
 all remaining conceptions of this article, always refering to the
 formula abbreviated by \Syntass.

 \begin{defin}\noem
     $\F_S = \bigcup\limits_{\gamma\in\SRT}\F_{\gamma}
     \qquad
     \sparenth{normal}{\vec{\sigma}\in\SRT^{\ell}} \quad
     \F_{\vec{\sigma}} = \Prod_{\indto i{\ell}} \F_{\sigma_i}
     $
 \end{defin}
%
%%####################################################################
%%# funl04.tex  |95-03-11|22:36|                                     #
%%####################################################################
 %         FUNL04.TEX
 %---------------------------+
 \section{Semantics of Functional Logic}
 %========================================%
 \begin{defin}\label{Strukt}\noem
     $\FnlSMG SM$
     \enspace:\enspace
     {$\frak{M}$ is a {$\Fnl$\em-type normal structure}\
     of {\em signature } $S$,
     iff the following conditions apply to it:
     }
 \end{defin}

 \begin{enumerate}^{\renewcommand{\labelenumi}{{\bf\arabic{enumi}}) }
         \renewcommand{\labelenumii}{\arabic{enumi}.\arabic{enumii}) }
         }
         %----------------------------------------------------------
     \item[{\bf i) }] $\fnlsigs=(\SRT,\SOP,\VSRT,\VAR,\sig)$
     \item[{\bf ii) }] $\M$ is a mapping defined on
         $\SRT \cup (\SRT \x \VSRT\stern) \cup \SOP$.
         This mapping assigns elements of $\SRT$ to corresponding ranges,
         members of $\SOP$ to symbol-interpretations
         (i.e. corresponding elements of or functions on such ranges
          or functionals in case of symbols that bind variables).
         To ordered pairs of $\SRT \x \VSRT\stern$ it assigns
         those components which determine the classes of functions
         admitted as arguments of the letter functionals.
     \item %=================================================(1)
         $\AQu{\alpha \in \SRT}\enspace
         \M(\alpha) \aliasgl \M_{\alpha} \neq \emptyset$
         and we automatically extend $\M$ to $\SRT\stern$:
         \begin{math}
             \AQu{\vec{\sigma}=
             \enangle{\sigma_i}_{\indto i\ell} \in \SRT\stern}
             \enspace
             \M(\vec{\sigma})\aliasgl\M_{\vec{\sigma}} \defgl
             \Prod_{\indto i\ell} \M_{\sigma_i}
         \end{math}
     \item %================================================(2)
             \begin{math}
             \AQu{\obl{\op}\in\SOP \komma \sigoop = \stdsig
             }\enspace
             (\text{using}\quad \M_{\alpha_i}^{\vec{\beta}_i} \aliasgl
                 \M(\alpha_i, \vec{\beta}_i)
             )
             \end{math}
             \\
             if $m=0$ then $\M\obl{\op}\in\M_{\gamma}$,
             \enspace otherwise
             \begin{math}
                 \Map \M\obl{\op}
                 :\Prod_{\indto im}\M_{\alpha_i}^{\vec{\beta}_i}
                 ->\M_{\gamma}.
             \end{math}\\
     \item %================================================(3)
         For arbitrary $\gamma\in\SRT \komma
         \vec{\sigma}=\enangle{\sigma_i}_{\indto i\ell} \in \VSRT^\ell$

         \begin{enumerate}^{\baselineskip=1.3\baselineskip}
             \item %------------------------------------------(3.1)
                 \begin{marray}[t]{rcccl}
                     \ell=0&\limp \M(\gamma,\vec{\sigma})
                         &= &\M(\gamma,\enangle{}) &= \M_{\gamma}
                     \\
                     \ell>0&\limp \M(\gamma,\vec{\sigma})
                     &\subseteq &\text{Map}(\M_{\vec{\sigma}},\M_{\gamma})%
                     &\defgl \M_{\gamma}^{\M_{\vec{\sigma}}}
                     \\
                 \end{marray}
                 \\
                 alias notation:\quad
                 $\M_\gamma^{\vec{\sigma}} \defgl \M(\gamma,\vec{\sigma})$

             \item[]{%% \large
             \underbar{\bf completion qualities of $\M_\gamma^{\vec{\sigma}}$}:
             }

             \item %------------------------------------------(3.2)
                 \begin{mArray}[t]{@{\extracolsep{3pt}}
                     rcccl@{\hspace{20mm}}r%
                     }%%%%%%%%%%%%%%%%%%%%%%%%%%%%%%%%%%%%%%%%%%%%%
                     \AQu{\reob{w} \in \M_\gamma}
                     &
                     \text{cst}^{\vec{\sigma}}_{\reob{w}}
                     &\defgl&
                     \funkd{\M\vec{\sigma}}{\M\gamma}{\revec{x}}{\reob{w}}
                     &\in\M_\gamma^{\vec{\sigma}}
                     &\quad \text{(constant funcs.)}
                     \\[3pt]
                     \AQu{\indto j\ell}
                     &
                     \text{pj}^{\vec{\sigma}}_j
                     &\defgl&
                     \funkd{\M\vec{\sigma}}{\M\sigma_j}%
                     {\enangle{\reob{x}_i}_i}{\reob x_j}
                     &\in\M_{\sigma_j}^{\vec{\sigma}}
                     &\text{(projections)}
                     \\
                 \end{mArray}
             \item %------------------------------------------(3.3)
             \begin{mArray}[t]{l}
                 \AQu{\vec{\varrho}\in\VSRT\stern}
                 \AQu{\revec{x}\in\M\vec{\sigma}}
                 \AQ{\reob{g}}
                 \\[3pt]
                 %\umklm(){%
                 \text{\underbar{if}}\quad
                 \reob{g}\in\M_\gamma^{\vec{\sigma}\verkett\vec{\varrho}}
                 \quad\text{\underbar{then}}\quad
                 %\boldsymbol{\limp}
                 \reob{g}_{\revec{x}} \defgl
                 \funkd{\M\vec{\varrho}\hspace{.6em}}{\M\gamma}%
                     {\revec{r}}{\;\reob{g}(\revec{x}\verkett\revec{r})}
                     \in \M_\gamma^{\vec{\varrho}}
                 %}
                 \hspace{16mm}\text{(partial fixing)}
             \end{mArray}
             \item %------------------------------------------(3.4)
             \begin{math}
                 \AQu{\obl{\op}\in\SOP \komma \sigoop = \stdsig}
                 \AQu{\reob g_1,...,\reob g_m}
             \end{math}
             \\
             \underbar{the premises}
             \qquad $m>0 \komma \ell>0$\qquad \underbar{and}
             \\
             \begin{math}
                 \AQu{\indto im}
                 \quad
             %    \begin{array}[t]{l}
                     \reob{g}_i \in
                     \M_{\alpha_i}^{\vec{\sigma}\verkett\vec{\beta}_i}
                     \quad
                     \text{ and introducing the auxiliary notation:}
                     \\
                     \reob{h}_i =
                     \begin{cases}
                         \reob{g}_i &\text{ if } r_i=0
                         \\
                     \funkd{\M\vec{\sigma}}%
                         {\text{Map}(\M\vec{\beta}_i\komma \M\alpha_i)}%
                         {\revec{y}}%
                         {{\reob{g}_i}_{\revec{y}} =
                         {\umklm[]{%
                             \revec{z} \mapsto
                             \reob{g}_i(\revec{y}\verkett\revec{z})}}
                         }
                         &\text{ if } r_i>0
                         \quad
                         \mbox{
             (3.3 implies $\reob{h}_i\revec{y}\in\M_{\alpha_i}^{\vec{\beta}_i}$)
                         }
                     \end{cases}
              %    \end{array}
             \end{math}
             \\
             \underline{imply}
             \\
             \begin{math}
                 \funkd{\M\vec{\sigma}}{\M\gamma}%

{\revec{y}}{\M\obl\op(\reob{h}_1(\revec{y}),\dots,\reob{h}_m(\revec{y}))}
                 \in
                 \M_{\gamma}^{\vec{\sigma}}
             \end{math}
             \hspace{30mm}(composition)
         \end{enumerate}

         \item[]%
     {\underline{\bf
         What $\M$ assigns to the fixed components of $S$}\enspace:}
         \item %==========================================(4)
         \begin{math}
             \M_{\prop} = \encurs{\M\obl{\falsum}\komma\M\obl{\verum}}
             = \mathbb{B} = \encurs{\AlgOp0, \AlgOp1}
         \end{math}
         and    \\
         \begin{math}
             \enangle{\M_{\prop}, \M\obl{\curlywedge}, \M\obl{\curlyvee},
                 \M\obl{\neg}, \M\obl{\wedge}, \M\obl{\vee}
             }
             = \AlgB =
             \enangle{\mathbb{B}, \AlgOp0, \AlgOp1,
                 \AlgOp{\com}, \AlgOp{\sqcap}, \AlgOp{\sqcup}}
         \end{math}

         forms a Boolean algebra with two elements,
         $\M\obl{\|\rightarrow\|}$ and $\M\obl{\|\leftrightarrow\|}$
         are represented by the (dependent) truth-operations
         $\sqimp_{\AlgB}$ und $\sqbipf_{\AlgB}$.
         $\M\obl{\GQ^{\alpha}}$ and $\M\obl{\PQ^{\alpha}}$
         are defined for $\alpha\in\VSRT$ as follows:
         \\[3pt]
         $%
             \M\obl{\GQ^{\alpha}}\komma \M\obl{\PQ^{\alpha}} \colon
             \M_{\prop}^{\enangle{\alpha}} \rightarrow \M_{\prop}
         $\quad
         for each $\theta \in \M_{\prop}^{\enangle{\alpha}}$ we stipulate
         \\
         if $(\boldsymbol{\GQ}\;\reob x \in \M_{\alpha})
              \enspace \theta\enangle{\reob x} = \AlgOp1
             $
             then
             $\M\obl{\GQ^{\alpha}}(\theta)=\AlgOp1$
             otherwise
             $\M\obl{\GQ^{\alpha}}(\theta)=\AlgOp0$;
         \\
         if $(\boldsymbol{\PQ}\;\reob x \in \M_{\alpha})
              \enspace \theta\enangle{\reob x} = \AlgOp1
             $
             then
             $\M\obl{\PQ^{\alpha}}(\theta)=\AlgOp1$
             otherwise
             $\M\obl{\PQ^{\alpha}}(\theta)=\AlgOp0$.
         \item %==========================================(5)
         \begin{math}
             \AQu{\alpha\in\SRT} \enspace
             \M\obl{\eqs{\alpha}} \colon
             \M_{\opaar{\alpha}{\alpha}} \rightarrow \M_{\prop}
             \komma\quad
             \opaar{\reob{x}}{\reob{y}} \mapsto
         %    \begin{cases}
                 \AlgOp1 \text{ if } \reob{x}=\reob{y}
                 \text{ or }
                 \AlgOp0 \text{ otherwise}
         %    \end{cases}
         \end{math}
 %%%        \\(also $\M\obl{\eqs{\alpha}} = \M_{\alpha}^2 \restr \Chi_{=}$)
 \end{enumerate}

 To extend a {\em structure } $\M$ into an {\em interpretation }
 of the language, i.e. to find an {\em evaluation } of {\em expressions }
 $\Expr_{\gamma}$ another approach than that based on
 {\em variables-assignments } as in {\em predicate logic } is required.
 The following definitions are prerequisites for the new approach.

 \begin{defin}\em\label{persp}
     $\Map\persp:\Expr_S->\PM(\VAR\stern).$
     \quad
 (Let $\fnlsigs \komma \Expr_S=\bigcup\limits_{\gamma\in\SRT} \Expr_{\gamma}$)
     \\
     For $\obl{e}\in\Expr_S$, $\persp \obl{e}$ denotes the set of all
     $\enangle{\obl{u_i}}_{\indto i\ell} \in \VAR\stern$
     such that all free variables of $\obl e$ are in
     $\clabst{\obl{u_i}}{1\le i \le \ell}$.
 \end{defin}

 \noindent
 We shall need
 a more technical approach in defining this conception using
 {\em syntactic induction}.
 If \Syntass\ (Def. \fullref{Sprache}) is assumed,
 then $\persp\obl{e}$ depends on $\persp\obl{a_i}$ as follows:
 \begin{displaymath}
 \begin{array}^{\Zzwi}% {\Zzwi\getMRCW{$m = 0$}}
     {>{\enspace}l|c||>{\enspace}l} % <-- Preambel
     \multicolumn{2}{l}{\text{{\bf cases}}} & \persp\obl{e} =
     \\ \hline
     %\multirow{2}{\MRCW}%
     {m = 0} &\oblt{op}\in\VAR
     &\clabst{\enangle{\obl{u_i}}_{\indto i\ell} \in \VAR\stern}{%
         \EQu{\indto j\ell}\; \obl{\op}=\obl{u_j}}
     \\
     \cline{2-3}
     &\oblt{op}\notin\VAR &\VAR\stern
     \\ \hline
     \multicolumn{1}{l}{\enspace m > 0}%
     &&\clabst{\vec{u}\in\VAR\stern}{%
             \AQu{\indto im}\;\vec{u}\verkett\vec{v_i}\in\persp\obl{a_i}}
             \vphantom{\Big\vert}
     \\ \hline
 \end{array}
 \end{displaymath}

 %% \newpagemittoc
 \begin{defin}\em\label{Lp}
     (Let $\fnlsigs  \komma \vec{u}\in\VAR\stern$)
 \end{defin}
 \begin{math}
     \parenth{normal}{\gamma\in\SRT} \enspace
         \Expr_{\gamma}[\vec{u}] =
         \clabst{\obl{e}\in\Expr_{\gamma}}{\vec{u}\in\persp\obl{e}}
     \qquad
     \parenth{normal}{\vec{\sigma}\in\SRT^\ell} \enspace
         \Expr_{\vec{\sigma}}[\vec{u}] = \textstyle
         \Prod_{\indto i\ell} \Expr_{\sigma_i}[\vec{u}]
 \end{math}

 \noindent
     $\Expr_{\gamma}[\vec{u}]$ is the set of expressions of
     $\Expr_{\gamma}$ whose free variables are among
     $\clabst{\obl{u_i}}{\indto i\ell}$ if
     $\vec{u} = \enangle{\obl{u_i}}_{\indto i\ell}$.
     \enspace
     $\Expr_{\gamma}[]$ therefore is the set of
      {\em closed} $\gamma-${\em expressions}.

 \begin{observ}\noem
     (for $\fnlsigs$): \qquad
     $\Expr_{\gamma} =
     \bigcup\limits_{\vec{u}\in\VAR\stern} \Expr_{\gamma}[\vec{u}]$
 \end{observ}
 \skipbasel{-1.3}

 \begin{observ}\noem\label{pGP}
     (for
     $\fnlsigs \quad \vec{u}=
      \enangle{\obl{u_i}}_{\indto i\ell} \in \VAR_{\vec{\sigma}}
      \quad \vec{\sigma}\in\VSRT^*
     $):
 \end{observ}
 \skipbasel{-0.9}
 \begin{displaymath}
     \enqt{e} \in \Expr_\gamma[\vec{u}] \leqv
     (\EQ{\oblt{op},m,\gamma,\vec{\alpha},\VVec{\beta},\vec{a},\VVec{v}})
     (\pSyntass)
 \end{displaymath}
 \noindent
     where \pSyntass\  (=perspective G.P.) can be obtained
     from \Syntass\ \apageref{Sprache} by modification of two conditions:
         if we change
         $\oblt{op} \in \SOP \cup \VAR$ into
         $\oblt{op} \in \SOP \cup
         \\ \cup
         \clabst{\obl{u_i}}{\indto i\ell}$
         \quad and \quad
         $\vec{a} \in \Expr_{\vec{\alpha}}$ into
         $\vec{a} \in
          \Prod_ {\indto im}\Expr_{\alpha_i}[\vec{u}\verkett\vec{v_i}]
         $\quad
          (each $\obl{a_i} \in \F_{\alpha_i}[\vec{u}\verkett\vec{v_i}]$).

 \begin{defin}\em
     \DKL{\underbar{Interpretation of the language into a structure}}{Interp}
     \\
 %\vspace{-1ex}
 Let
 %\begin{minipage}[t]{.9\textwidth}
     $\FnlSM SM \komma \gamma\in\SRT \komma
      \obl{e}\in\Expr_{\gamma} $
      and \Syntass\ be assumed.
 %\end{minipage}\\[1ex]
 The evaluation $\obl{e}_{\M}$ of $\obl{e}$ is defined to be
 a function on $\persp\obl{e}$. Let
 $%\begin{math}
     \vec{u}=\enangle{\obl{u_i}}_{\indto i\ell} \in \persp\,\obl{e}
     \komma
     \AQu{\indto i\ell} \usig\obl{u_i}=\obl{\sigma_i}
     \komma
     \vec{\sigma}=\enangle{\sigma_i}_{\indto im}
 $. %\end{math}.
 Then $\obl{e}_{\M}(\vec{u})$ is defined inductively:
 \end{defin}
 \skipbasel{-1.2}
 \begin{displaymath}
 %\hspace*{-3em}
 \begin{array}^{\Zzwi \newcommand{\mzw}[1]{\multicolumn{2}{c||}{#1}}}%
 {c|c|c||>{\enspace}l}
     \multicolumn{3}{c}{\text{{\bf cases}}} & \obl{e}_{\M}(\vec{u}) =
     \\ \hline
     \ell=0 &\mzw{m=0} &=\M\oblt{op}
     \\
     \cline{2-4}
     \begin{sizemath}{\small}
         (\vec{u}=\enangle{})
     \end{sizemath}
     &\mzw{m>0}
     &=\M\oblt{op}
       \parenth{big}{\enangle{\obl{a_i}_{\M}(\vec{v_i})}_{\indto im}}
     \\ \hline
     \ell>0 &m=0 &\oblt{op} \in \VAR
     &=%
         \text{pj}_k^{\vec{\sigma}} =
         \funkd{\M\vec{\sigma}}{\M\sigma_k}{\enangle{\reob{x}_i}_i}{\reob{x}_k}
         \quad
      \begin{Array}^{\footnotesize}{l}
         \text{where } k=
         \\
         \mathop{\operatorname{max}j}\limits_{\indto j\ell}(\obl{u_j}=\oblt{op})
      \end{Array}
     \\
     \cline{3-4}&&&
     \\[-2ex]
     &&\oblt{op} \in \SOP &=\text{cst}^{\vec{\sigma}}_{\M\oblt{op}} =
         \funkd{\M\vec{\sigma}}{\M\gamma}{\revec{x}}{\M\oblt{op}}
     \\[3ex]
     \cline{2-4}
     &\mzw{}&
     \\[-2ex]
     &\mzw{m>0}&=%
         \funkd{\M_{\vec{\sigma}}}{\M_{\gamma}}%
         {\revec{x}}{\M\oblt{op}(\enangle{\reob{h}_i(\revec{x})}_{\indto im})}
         \quad
         \parbox{18mm}{\footnotesize
             with $\reob{h}_i$ defined below by *) }
     \\
     \hline
 \end{array}
 \end{displaymath}
 \begin{math}
     \text{\small *)}\enspace
     \begin{Array}[t]{l}
         (\text{case }\ell>0,m>0)\AQu{\indto im}
         \reob{h}_i\colon\M_{\vec{\sigma}}\to\M_{\alpha_i},
         \enspace
         \text{if } r_i = 0 \colon
         \reob{h}_i\colon\revec{x} \mapsto \obl{a_i}_{\M}(\vec{u})(\revec{x})
          \\
         \text{if } r_i > 0 \text{ then }
         \reob{h}_i\colon\revec{x} \mapsto
         \parenth{big}{\obl{a_i}_{\M}(\vec{u}\vec{v_i})}_{\revec{x}} =
         \funkd{\M_{\vec{\beta_i}}}{\M_{\alpha_i}}%
         {\revec{y}}{\obl{a_i}_{\M}(\vec{u}\vec{v_i})(\revec{x}\revec{y})}.
     \end{Array}
 \end{math}

 \begin{eigen}\noem \quad
     \begin{math}
         \obl e \in \Expr_{\gamma}[\vec{u}] \land
         \vec{\sigma} \in \VSRT\stern
         \land
         \vec{u} \in \VAR_{\vec{\sigma}}
         \limp
         \obl{e}_{\M}(\vec{u}) \in \M^{\vec{\sigma}}_{\gamma}
     \end{math}
 \end{eigen}

 \begin{PROOF}{(+ Remark)}
     This \eigensch\ is already required for the argument expressions
     $\obl{a_i}$ of the preceding definition (\ref{Interp})
     to assert that $\enangle{\reob{h}_i(\revec{x})}_{\indto im}$
     belongs to the domain of $\M\oblt{op}$
     (this assertion also requires (3.3) of \ref{Strukt} {\em def.}).
     Conditions \ref{Strukt}(3) imply that the above
     {\it\eigensch } propagates from the $\obl{a_i}$ to $\obl e$;
     so {\em syntactic induction } ensures its validity and
     any circularity of \ref{Interp} {\em def.}
     that might result from presupposing it (for $a_i$) is avoided
     as well.
 \end{PROOF}

 \begingroup %>>>>>>>>>>>>>>>>
 \renewcommand{\N}{\frak{N}}

 \begin{observ}\em\label{eq.eval}
     If $\M,\N \in \fnlstruct \spa
     \AQu{\gamma\in\SRT} \M_{\gamma} = \N_{\gamma}
     $\quad
     and
     \\
     $\AQu{\oblt{op}\in\SOP \kom \sigoop=\stdsig}
      \AQu{\revec{h}\in
             \Prodim\M^{\vec{\beta_i}}_{\alpha_i}
             \cap
             \Prodim\N^{\vec{\beta_i}}_{\alpha_i}
         }
       \M\oblt{op}(\revec{h})=\N\oblt{op}(\revec{h})
     $
     \\
     (for $m=0,\revec{h}=\emptyseq:\quad
       \M\oblt{op}(\revec{h})=\N\oblt{op}(\revec{h})$)
     \hfill then
     $\AQu{\obl{e}\in\bigcup\limits_{\gamma\in\SRT}\Expr_{\gamma}}
         \obl e_{\M} = \obl e_{\N}
     $
 \end{observ}

 \begin{proof} syntactic induction on $\obl{e}$
 \end{proof}

 \begin{thms}{Conclusion}\em\label{FullStr}
     If $\overline{\M}$ is characterized by
     $\AQ{\gamma}\;\BM_{\gamma}=\M_{\gamma}
     \quad
     \AQ{\gamma,\vec{\sigma}}\;
     \BM_{\gamma}^{\vec{\sigma}}={\BM_{\gamma}}^{\BM_{\vec{\sigma}}}
     $
     \\
     and\quad
     $\AQ{\oblt{op}}\,\parenth{big}{\Prodim\M_{\alpha_i}^{\vec{\beta_i}}}
         \restr \BM\oblt{op}
         = \M\oblt{op}
     $
     \quad then\quad
     $\AQu{\obl{e}\in\bigcup\limits_{\gamma\in\SRT}\Expr_{\gamma}}
         \obl{e}_{\BM}=\obl{e}_{\M}
     $.
 \end{thms}
 \endgroup %%%<<<<<<<<<<<<<<<<

%% file: funl5-6.tex
% Logic Eprints
%Submitted 0644 Tue Mar 14, 1995 by: a8121dab@helios.edvz.univie.ac.at (josef schoenbrunner )
%logic/schoenbrunner/fun-logic-completeness
%logic/schoenbrunner/fun-logic-completeness/funl5-6.tex
%

% FUNL5-6.TEX = FUNL05.TEX + FUNL06.TEX $1
%%####################################################################
%%# funl05.tex  |95-03-12|10:43|                                     #
%%####################################################################
 % FUNL05.TEX %
 % ~~~~~~~~~~~%
 \section{Syntactic Matters and the Calculus of Functional Logic}
 The {\em logical axioms } depend on the syntactic notions
 {\em free variables, bound variables } of an expression
 $\obl{e}\in\Expr_{\gamma}$, {\em substitutability } and {\em substitution}\
 (of a variable in an expression for some term).

 \begin{thms}{Notation}\noem\label{vcmps}
    $\vcmps a$ symbolizes the set of components of an arbitrary finite
    sequence $\vec{a}$
    (if $\vec{a}=\enangle{\obl{a_i}}_{\indto i\ell}$,
    then $\vcmps a = \clabst{\obl{a_i}}{\indto i\ell}$)
 \end{thms}

 \begin{defin}\label{fbV}
    \em$\fv\obl{e}\komma\gv\obl{e}$
    (free and bound variables in $\obl{e}\in\Expr_{\gamma}$).
    Provided that $\obl{e}\in\Expr_{\gamma}$ and \Syntass\
    we define inductively:
 \end{defin}
 \vspace{-1\baselineskip} % {-1ex}
    \begin{displaymath}\textstyle
    %\hspace*{-2.5em}
    \begin{Array}^{\Zzwi}{@{\extracolsep{1em}}c|c|c}
       \multicolumn{1}{c}{\text{cases}}
       &\multicolumn{1}{c}{\fv\obl{e}=}
       &\multicolumn{1}{c}{\gv\obl{e}=}
       %\text{cases}&\fv\obl{e}=&\gv\obl{e}=
       \\ \hline %---------------------------------
       m=0\enspace &\encurs{\oblt{op}} \cap \VAR &\emptyset
       \\
       m>0\enspace &%
       \bigcup\limits_{\indto im}
       (\fv\obl{a_i} \setminus \vcmps{v_i}
       %%% \clabst{\obl{v_{ij}}}{\indto{j}{r_i}}
       )
       &%
       \bigcup\limits_{\indto im}
       (\gv\obl{a_i} \cup \vcmps{v_i}
       %%% \clabst{\obl{v_{ij}}}{\indto{j}{r_i}}
       )
       \\ \hline %---------------------------------
    \end{Array}
    \end{displaymath}

 \underbar{Remark}:
    If $\oblt{op} \in \VAR$ then
    $\encurs{\oblt{op}} \cap \VAR = \encurs{\oblt{op}}$,
    otherwise
    $\encurs{\oblt{op}} \cap \VAR = \emptyset$.

 \begin{bemerk}\em\label{frV+pers}
    $\begin{array}[t]{L@{\enspace}rl}
       (1)
       &\obl{e} \in \Expr_{\gamma}[\vec{u}]
       &\leqv
       \obl{e} \in \Expr_{\gamma} \land
       \fv\,\obl{e} \subseteq \vcmps u
       %%% \clabst{\obl{u_p}}{\indto p{\ell_{\vec{u}}}}
       \\(2)
       &\obl{e} \in \Expr_{\gamma}[\vec{u}]
       &\limp
       \parenth{big}{%
          \vcmps u
          \cap W = \emptyset
          \limp
          \fv\,\obl{e}
          \cap W = \emptyset
       }
     \end{array}
    $
 \end{bemerk}

 \begin{defin}\em\label{substabl}
    \underbar{Substitutability}:
    \enspace
    \begin{math}
       \text{Subb} \subseteq \Expr_{S}\x\VAR\x\Expr_{S}
    \end{math}
 \end{defin}
    \skipbasel{-1}
    \begin{multline*}
       \Syntass \limp
       \sbbobl{d}{x}{e} \leqv
       \\
       \leqv \AQu{\indto im}\parenth{Big}{%
          \obl{x}\in
          \vcmps{v_i} \lor
          \parenth{big}{\sbbobl{d}{x}{a_i} \land
             \fv\obl{d} \cap
             \vcmps{v_i} = \emptyset}}
    \end{multline*}

 \begin{thms}{Notation}\noem
    ${\obl{\sub exd} \aliasgl \obl{e^x_d}}$
    \enspace
    denotes the result of replacing
    each {\em free occurring } $\obl x$ by $\obl d$ applied to $\obl e$.
    This is a special case of next Definition (with $l=1$).
 \end{thms}

 \begin{defin}\noem\label{subst}
       $\obl{\sub{e}{\vec{x}}{\vec{d}}}
        \aliasgl \obl{e^{\vec{x}}_{\vec{d}}}$
    \enspace (for
    \begin{math}
       \obl e\in\Expr_{\gamma}
       \komma
       \vec{\sigma} = \enangle{\sigma_i}_{\indto i{\ell}}
       \in\VSRT^{\ell}
       \komma
       \vec{x} = \enangle{\obl{x_i}}_{\indto i{\ell}}
          \in %%%%% \Prod_{\indto i{\ell}}\VAR_{\sigma_i} \aliasgl
          \VAR_{\vec{\sigma}}
       \komma
       \vec{d} = \enangle{\obl{d_i}}_{\indto i{\ell}}
          \in \Expr_{\vec{\sigma}}
    \end{math}
    )
    denotes the result of simultaneously replacing each {\em free occurring }
    $\obl{x_i}$ by $\obl{d_i}\quad(\indto i{\ell})$  applied to $\obl e$.
    If a variable appears more than once within the sequence $\vec{x}$,
    the rightmost $d_i$ of the corresponding position replaces the variable.
    This is defined inductively: if \Syntass\  is supposed, then
 \end{defin}
 \skipbasel{-1.5}
 \begin{displaymath}%\hspace*{-2.4em}
    \begin{array}^{\Zzwi\newcommand{\mzw}[1]{\multicolumn{2}{c|}{#1}}}{c|c|l}
       \multicolumn{2}{c}{\text{{\bf cases}}}
       &\obl{\sub{e}{\vec x}{\vec d}} =
       \\ \hline %---------------------
       m=0 &\oblt{op} \in \vcmps x
       &=\obl{d_k}\enspace
       \mbox{{\small with
          $k=\operatorname{max}_{\indto j{\ell}}(\oblt{op}=\obl{x_j})$
          }}
       \\
       \cline{2-3}%--------------------
       &\text{otherwise} &=\oblt{op}\qquad(= \obl e)
       \\ \hline %---------------------
       \mzw{m>0}
       &=\obl{\opexx{v}{a}{\psubvv{x}{v_i}{d}{v_i}}
       }
       \\ \hline %---------------------
    \end{array}
 \end{displaymath}

 \noindent
    Substitution $[\vec x \verkett \vec v_i \zu \vec d \verkett \vec v_i]$
    differs from $[\vec x \zu \vec d]$ exactly if
    the sequences $\vec x$ and $\vec v_i$ have common members.
    If $\obl{x_i}=\obl{v_{ij}}$ then the replacement of
    $x_i$ by $d_i$ is prohibited in the substitution
    $[\vec x \verkett \vec v_i \zu \vec d \verkett \vec v_i]$.
    This prevents replacing bound variables.

 \begin{konv}\label{id.in.sub}
    Inside of $[...\gets...]$ an identifier $a$ for any expression
    $\obl{a} \in \Expr_{\alpha}$ denotes the sequence $\enangle{\obl{a}}$
    of length 1. E.g. $e\psub{\vec{x}y\vec{z}}{\vec{a}b\vec{c}}$ is to be
    read as
    $e\psub{\vec{x}\enangle{\obl{y}}\vec{z}}{\vec{a}\enangle{\obl{b}}\vec{c}}$.
 \end{konv}

 \begingroup %>>>>>>>>>>>>>>>>
    \newcommand{\If}{$. If\enspace $}
    \newcommand{\nlIf}{\SPLITMBOX{.\\ If\quad }}
    \newcommand{\Then}{\SPLITMBOX{\spa then }}
    \def\PSTEXT{\hfill}
    \newcommand{\NZfL}{\SPLITMBOX{\\}}
    \newenvironment{LEMMA}[2][0]{\SETHBOX[#1]
       \lemma\noem\label{#2}
          \enspace Let\enspace
          \BEGINHBOX$%
    }{$\ENDBOX.\PSTEXT\endlemma}
    \newcommand{\VVECTYPE}[2]{\vec{#1}\in\VAR_{\vec{#2}}}
    \newcommand{\EVECTYPE}[2]{\vec{#1}\in\Expr_{\vec{#2}}}
 %==============================================
 \vspace{1ex}
 \noindent
 Laws of \textit{substitution} enumerated within
 the subsequent five \textit{lemmas}
 shall prove to be essential prerequisites for propositions
 concerning the {\em term structure } obtained from a consistent {\em theory}
 in our final section.
 The proofs of these (intuitively clear) {\em lemmas } \ref{HS1} to \ref{HS5}
 mainly rely on {\em syntactic induction } using \GENP.

 \begin{LEMMA}^{\def\PSTEXT{\hfill This is a special case of the next}}{HS1}
    \obl{e} \in \Expr_{\gamma} \quad
    \vec{\varrho},\vec{\sigma} \in \VSRT^* \quad
    \vec{u} \in \VAR_{\vec{\sigma}} \quad
    \vec{c} \in \Expr_{\vec{\sigma}}
    \nlIf
    \fv\,\obl{e} \cap \invec{u} = \emptyset
    \Then
    \obl{e\psubv uc} = \obl{e}
 \end{LEMMA}

 \begin{LEMMA}[*]{HS2}
    \obl{e} \in \Expr_{\gamma}%   %[\vec{x}\vec{y}]
    \quad
    \VVECTYPE y\eta \quad \VVECTYPE u\varrho
    \quad
    \vec{\eta},\vec{\varrho} \in \VSRT^*
    \quad
    \EVECTYPE r\varrho \quad \EVECTYPE s\eta.
    \SPLITMBOX{\Hfill If\quad}
    \fv\obl{e} \cap \vcmps u  \subseteq \vcmps y
    \quad \text{then} \quad\hfill
    \obl{e\psubvv uyrs} = \obl{e\psubv ys}
 \end{LEMMA}
 \begin{proof}
    From the premises (i)
    $\fv\obl{e} \cap \vcmps u  \subseteq \vcmps y$
    \quad (ii) \GENP and (iii) \textit{induction hypotheses}\quad
    we shall infer the succedent \ileft=\iright.
   \begin{liste}[\labelwidth=1.7em]{\ItemCase}%%{\leftmargin=0em}
      \item $m=0$.\qquad Then
    $\obl{e}=\oblt{op}$ and $\fv\obl{e}=\encurs{\obl{e}} \cap \VAR$.
    By (i) and the fact, that $\vcmps{u} \subseteq \VAR$
    we obtain (iv)
    \begin{math} % \begin{equation}\tag{iv}
        \encurs{\obl{e}} \cap \vcmps{u} \subseteq \vcmps{y}
    \end{math}   % \end{equation}
    \begin{liste}{{\CaseLogo}}
       \item   $\obl{e} \in \cmpset{\vec{u}\vec{y}}$.
          \quad Then (iv) implies (v) $\obl{e} \in \vcmps y$.
          According to
          \Vref{subst}{def., case \em $\oblt{op}\in\cmpset{...}$,}
          we have
          \begin{monolist}{(vi)}
             $\ileft=\obl{e\psubvv uyrs}=\pj_k(\vec{r}\vec{s})$,
             where $k$ is the maximal within range
             ${1 \dots \ell_{\vec{\varrho}\vec{\eta}}}$, so that
             $\pj_k(\vec{u}\vec{y})=\obl e$. \hfill (v) implies that
            \item[(vii)] $k=\ell_{\vec{\varrho}}+j$,
             where j is the maximal within range
             ${1 \dots \ell_{\vec{\eta}}}$ so that
             $\pj_j(\vec{y})=\obl{y_j}=\obl{e}$.
          \end{monolist}
          (vi)+(vii) yield
          $\ileft=\pj_{\ell_{\vec{\varrho}}+j}(\vec{r}\vec{s})=\obl{s_j}$
          and
          $\iright=\obl{e\psubv ys}=\obl{s_j}$; $\ileft=\iright$.
       \item   $\obl{e} \notin \cmpset{\vec{u}\vec{y}}$.
          According to
          \Vref{subst}{def., case \em $\oblt{op}\notin\cmpset{...}$,}
          $\ileft=\iright$.
    \end{liste}
    \item $m \neq 0$.
    \begin{liste}[\usecounter{enumi}\setcounter{enumi}{3}]{(\roman{enumi})}
       \item
       \begin{math}
          (\fv\obl{a_i} \setminus \vcmps{v_i}) \cap \vec{u}
          \subseteq \vec{y}
       \end{math}
       \hfill
       from (i), \Vref{fbV}{def.of \em $\fv$}
       \item
       \begin{math}
          \fv\obl{a_i} \cap \vcmps{u} \subseteq
          \vcmps{y} \cup \vcmps{v_i} = \cmpset{\vec{y}\vec{v_i}}
       \end{math}
       \hfill as $A \subseteq B$ implies
       $A \cup \vcmps{v_i} \subseteq B \cup \vcmps{v_i}$
       \item
       \begin{math}
          \ileft = \obl{e\psubvv uyrs}
          \begin{Array}[t]{cl}
             \undertext{\ref{subst}}=&%
             \obl{\opexb[\psub{\vvvec uy{v_i}}{\vvvec rs{v_i}}]va}
             \\ \undertext{(iii)+(v)}=&%
             \obl{\opexb[\psubvv y{v_i}s{v_i}]va}
             \\ \undertext{\ref{subst}}=&%
             \obl{e\psubv ys} = \iright
          \end{Array}
       \end{math}
    \end{liste}
   \end{liste}
   \skipbasel{-2}
 \end{proof}

 \begin{LEMMA}[*]{HS3}
    \obl{e} \in \Expr_{\gamma} \quad
    \vec{x} \in \VAR_{\vec{\xi}} \quad
    \vec{y} \in \VAR_{\vec{\eta}} \quad
    \vec{\xi},\vec{\eta} \in \VSRT^* \quad
       \vec{r} \in \Expr_{\vec{\xi}} \quad
       \vec{s} \in \Expr_{\vec{\eta}}
          \SPLITMBOX{\\ \Hfill If\quad}
       \UnVec{\fv}{r} \cap (\invec{y} \cup \gv\obl{e}) = \emptyset
    \Then
    \obl{e\psubvv xyrs} = \obl{e\psubvv xyry \psubv ys}
 \end{LEMMA}
 \noindent\textbf{Remark }
 If $\vec{r} \in \Expr_{\vec{\xi}}[]$, then the last condition is true.
 \begin{proof}
    From the premises (i,ii,iii) we shall infer $\ileft=\iright$
    (succedent of the lemma). (i) \textit{premises of lemma}; \quad
    (ii) \GENP; \quad (iii) \textit{induction hypothesis}
   \begin{liste}[\labelwidth=1.7em]{\ItemCase}%%{\leftmargin=0em}
      \item $m=0$.\qquad Then $\obl{e}=\oblt{op}$
    \begin{liste}{\CaseLogo}
       \item $\obl{e} \notin \cmpset{\vec{x}\vec{y}}$.
          Then $\ileft=\obl{e}=\iright$ (according to \Vref{subst}{def.})
       \item $\obl{e} \in \vcmps{x}\setminus\vcmps{y}$.
          Then $\ileft=\obl{r_m}$, where $m$ is the maximal
          $m \le \ell_{\xi}$ so that $\obl{e}=\obl{x_m}$.
          Hence $\iright=\obl{r_m\psubv ys}$.
          (i) and \ref{HS1} yield $\iright=\obl{r_m}$;
          $\ileft=\iright$
       \item $\obl{e} \in \vcmps{y}$.
          Then $\ileft=\obl{s_m}$ where $m$ is the maximal so that
          $\obl{e}=\obl{y_m}$, hence
          $\iright=\obl{y_m\psubv ys}=\obl{s_m}$, \ileft=\iright.
    \end{liste}
    \item $m>0$.
       Now the premises (ii),(iii) become relevant.
       Again we infer \ileft=\iright.
       \begin{liste}[\setcounter{enumi}{4}]{(\roman{enumi})}
          \item[(iv)]
          \begin{math}
             \ileft=\obl{\opexb[\psub{\vvvec xy{v_i}}{\vvvec rs{v_i}}]va}
          \end{math}
          \\
          As $\vcmps{v_i} \subseteq \gv\obl{e}$ (i) implies
          $\UnVec{\fv}{r} \cap
           (\cmpset{\vvec{y}{v_i}} \cup \gv{e})=\emptyset$,
           then (iii) yields
          \item[(v)]
          \eqobl{a_i\psub{\vvvec xy{v_i}}{\vvvec rs{v_i}}}%
          {a_i\psub{\vvvec xy{v_i}}{\vvvec ry{v_i}}\psubvv{y}{v_i}{s}{v_i}}
          \\
          \Hfill Substituting in (iv) and application of \ref{subst} yield
          \item[(vi)]
          \begin{math}
             \ileft
            \begin{Array}[t]{cl}
             =&%
             \obl{\opexb[\psub{\vvvec xy{v_i}}{\vvvec ry{v_i}}]va\psubv ys}
             \\
             \undertext{\ref{subst}}=
             &\obl{e\psubvv xyry \psubv ys}=\iright
            \end{Array}
          \end{math}
       \end{liste}
   \end{liste}
   \skipbasel{-1}
 \end{proof}

 \begin{LEMMA}[*]{HS3a}
    \obl{e} \in \Expr_{\gamma} \quad
    x \in \VAR_{\varrho} \quad
    \vec{y} \in \VAR_{\vec{\eta}} \quad
    \varrho,\vec{\eta} \in \VSRT^* \quad
    \obl{r} \in \Expr_{\varrho}
 \end{LEMMA}
 \EquaNull \vspace{-1ex} % \skipbasel{-1}
 \hfill
 \begin{math}
    \obl{e\psub{x\vec{y}}{r\vec{y}}}
    =
    \begin{cases}
       \obl{e}
       &\text{if } \obl{x} \in \invec{y}
       \\
       \obl{e\psub{x}{r}} \quad
       &\text{if } \obl{x} \notin \invec{y}
    \end{cases}
 \end{math}
 \hspace{1.3in}
 \begin{proof}
    Again rely on \GENP and use induction,
    the result then comes immediately from \ref{subst}\ \textit{def.}
 \end{proof}

 For the purpose of proving the next lemma we observe
 \begin{sublemma}\noem
    Let
    \begin{math}
       \obl{e} \in \Expr_{\gamma} \quad
       \vec{x} \in \VAR_{\vec{\xi}} \quad
       \vec{y} \in \VAR_{\vec{\eta}} \quad
       \vec{\xi},\vec{\eta} \in \VSRT^* \quad
       \vec{r} \in \Expr_{\vec{\xi}}.
       \\ \Hfill\text{If }
       \disjvecs xy \quad\text{then}\quad
       \obl{e\psubvv xyry} = \obl{e\psubv xr}
    \end{math}
 \end{sublemma}

 \begin{LEMMA}[*]{HS5}%
    \obl{e} \in \Expr_{\gamma} \quad
    \vec{x} \in \VAR_{\vec{\xi}} \quad
    \vec{y} \in \VAR_{\vec{\eta}} \quad
    \vec{\xi},\vec{\eta} \in \VSRT^* \quad
    \vec{c} \in \Expr_{\vec{\xi}}[] \quad
    \vec{d} \in \Expr_{\vec{\eta}}[]
    \nlIf
    \disjvecs xy
    \quad\text{then}\quad
    \obl{e\psubvv xycd} = \obl{e\psubv xc \psubv yd}
    = \obl{e \psubv yd \psubv xc}
 \end{LEMMA}

 \begin{proof}
    Proving the \textit{sublemma} is as easy as for the previous lemma.
  %   the same remark as for the preceding lemma applies.
    For the current \textit{lemma} two equations are considered.
    From \ref{HS3}+\textit{remark} and the preceding \textit{sublemma} we
    immediately obtain the first
    $(\obl{e\psubvv xycd} = \obl{e\psubv xc \psubv yd})$.
    in order to enable the induction step for the second equation
    we prove more generally for arbitrary $\vec{u}\in\VAR^*$
    \begin{equation}\tag{*}
    %   \eqobl{e\psubvv xucu \psubvv yudu}{e\psubvv yudu \psubvv xucu}
       \eqobl{e\psub{\vvvec xyu}{\vvvec cdu}}{e\psub{\vvvec yxu}{\vvvec dcu}}
    \end{equation}
    Assume \textit{premises of the lemma},\quad
    \GENP \quad and \textit{induction hypothesis}.
   \begin{liste}[\labelwidth=1.7em]{\ItemCase}%%{\leftmargin=0em}
      \item $m=0$.\qquad Then $\obl{e}=\oblt{op}$
    \begin{liste}{\CaseLogo}
       \item $\obl{e} \notin \cmpset{\vvvec xyu}$.
          Then $\ileft=\obl{e}=\iright$
          (according to \Vref{subst}{def.}, 2nd line of the table)
       \item $\obl{e} \in \vcmps{u}$.
          Then $\eqobl{e}{u_p}$.
          From \ref{subst} \textit{def.} (1st line of the table)
          and the observation,
          that the rightmost occurrence of $\obl e$ within $\vvvec xyu$
          belongs to $\vec{u}$ we infer
          $\ileft=\pj_{\ell_x+\ell_y+p}(\vvvec xyu)=\obl{u_p}=\obl{e}$.
          and by obvious symmetric consideration $\iright=\obl{e}$.
       \item $\obl{e} \in \vcmps{x} \setminus \vcmps{u}$.
          Then, as $\disjvecs xy$ is assumed,
          $\obl{e} \notin \cmpset{\vvec yu}$ and
          \ref{subst} \textit{def.}, 1st line of table
          yields
          $\ileft=\obl{c_m}$, where $m$ is the maximal
          $m \le \ell_{\xi}$ so that $\obl{e}=\obl{x_m}$.
          Then  $\ell_{\vec{y}}+m$ is the corresponding position
          for \iright, i.e.
          \begin{math}
             \iright=\pj_{\ell_{\vec{y}}+m}(\vvvec dcu)
             =\pj_m(\vec{c})=\obl{c_m}
          \end{math};
          \ileft=\iright.
       \item $\obl{e} \in \vcmps{y}\setminus\vcmps{u}$
          is obviously similar to the preceding case
          (exchange $x$ with $y$ as well as \ileft with \iright)
    \end{liste}
    \item $m>0$.
       From \GENP and
       \ref{subst} \textit{def.}, 3rd line of table
       we obtain
          \begin{mArray}{rl}
             \ileft &=
             \obl{\opexb[%
                \psub{\vvec xy \vvec{u}{v_i}}{\vvec cd \vvec{u}{v_i}}
                ]va}
             \\
             \iright &=
             \obl{\opexb[%
                \psub{\vvec yx \vvec{u}{v_i}}{\vvec dc \vvec{u}{v_i}}
                ]va}
          \end{mArray}.
       \\
       Replacing the argument-expressions according to the
       \textit{induction hypothesis}
       ($\obl{a_i}$ in place of $\obl{e}$
         and $\vvec{u}{v_i}$ in place of $\vec{u}$) yields \ileft=\iright.
   \end{liste}
   \skipbasel{-1}
 \end{proof}
 \endgroup %%%<<<<<<<<<<<<<<<<

 \begin{mkonv}\noem
    (1) We write $\obl{\AQ{x}}$ for $\obl{\GQ^{\xi}x}$
    if $\obl{x} \in \VAR_{\xi}$,
    and (2) $\obl{\AQ{\vec{z}}}$ for $\obl{\AQ{z_1}\dots\AQ{z_{\ell}}}$
    if $\vec{z}=\enangle{\obl{z_i}}_{\indto i{\ell}}$.
    This includes {\em case } $\ell=0$, as for $\vec{u}=\emptyseq$
    \quad
    $\obl{\AQ{\vec{u}}\varphi}=\obl{\varphi}$ is stipulated.
    (3) $\obl{a=b}$ is an alias for $\obl{a \eqs{\gamma} b}$
    if $a,b \in \Expr_{\gamma}$.
 \end{mkonv}

 \begin{thms}{Calculus of Functional Logic}\label{CAL}
    \noem \quad $\KF=\KF_S$ is defined for $\fnlsigs$
    as a triple of the component sets
    {\em formulae, axioms } and {\em rules } specified as follows:
 \end{thms}
 \begingroup %***************************
    \newcommand{\component}[1]{\textbf{\textit{\underbar{#1}}}\xspace}
    \newcommand{\subcompon}[1]{{\it\underbar{#1}}\enspace}
    \newcommand{\mobl}[1]{$\obl{#1}$}
    \let\ssep\komma % \newcommand{\ssep}{\hspace{.2em};\enspace}
    \jnmusk=\thinmuskip
 %*********************************************************
    \component{Formulae}
       $\Expr_{\prop}$ \quad (s. \fullref{Sprache}).\quad
    \component{Axioms}
       \subcompon{Propositional Tautologies};%
       \\%\quad
       \subcompon{Predicate Logic Axioms}
          %-----------------------------%
          Like in {$1^{\text{st}}$ \em order predicate logic }
          but related to the extended notions of
          {\em free vars., bound vars., substitutability}
          and {\em substitution}.
          The Axioms are
             \mobl{\AQ{x}\,\phi \limp \phi^x_a}
             \ssep
             \mobl{\phi^x_a \limp \EQ{x}\,\phi}
             \enspace ({\em meta condition } $\sbbobl{a}x{\phi}$ provided)
             \\
             \mbox{%
             \mobl{\AQ{x}(\psi \limp \phi)
                \limp (\psi \limp \AQ{x}\,\phi)},
             \enspace
             \mobl{\AQ{x}(\phi \limp \psi)
                \limp (\EQ{x}\,\phi \limp \psi)}
             \\[2pt] \Hfill
             (if $\obl x \notin \fv\obl{\psi}$).%
             }
       \\
       \subcompon{Equality Axioms}\label{eqax}
          \\
             (I1) \mobl{a \eqs{\alpha} a}
             \hfill
             $%------------------
             (\text{for }
             \alpha \in \SRT
                \komma
                \obl{a} \in \Expr_{\alpha})
             $%------------------
             \\
             (I2) \mobl{%
                \AQvecz
                \enspace
                \sub{b_1}{\vec{x}}{\vec{z}}
                \eqs{\alpha_i}
                \sub{b_2}{\vec{y}}{\vec{z}}
                \;
                \boldsymbol{\rightarrow}
                \;
                \OPEXDG{x}{b_1} \eqs{\gamma} \OPEXDG{y}{b_2}
                }
    \\
    \Hfill
      \begin{minipage}{.96\textwidth}
   %   \begingroup %>>>>>>>>>>>>>>>>
       \newcommand{\sep}{\hspace{1em}}
       where
       \begin{math}
          \oblt{op} \in \SOP \komma
          \sigoop = \sigopx \sep
          \vec{x},\vec{y},\vec{z} \in
             \VAR_{\vec{\beta_i}}
          \sep
          b_1,b_2 \in \Expr_{\alpha_i}
          \\
          \vcmps z
          \cap (\fv\obl{b_1} \cup \fv\obl{b_2})
          = \emptyset
       \end{math}
   %   \endgroup %%%<<<<<<<<<<<<<<<<
       \quad and $\Delta$ and $\Gamma$
       may be further arguments or empty.
      \end{minipage}
    \\[1ex]
    \component{Rules}
          \subcompon{Detachment}
          $\obl{\phi}\komma \obl{(\phi \rightarrow \psi)}
             \mapsto \obl{\psi}
          $
          and \subcompon{Generalization}:
          $\obl{\phi} \mapsto \obl{\AQ{x}\;\phi}$
          \ifopts0 \\ \hspace*{30mm} \else \quad \fi
          {\small $%
          (\;\obl{\phi},\obl{\psi} \in \Expr_{\prop}
           \komma \obl{x} \in \VAR
           \;)$}

 \underbar{Remark}:
 For $\alpha \in \VSRT$ scheme (I1) could be replaced by a single axiom
    \mobl{x \eqs{\alpha} x} (for one fixed $\obl x \in \VAR_{\alpha}$).
 \endgroup %***************************

 \begin{thms}{Notation}\em
    $\obltup{\varphi}m \ded{S} \obl{\psi} \quad (m\ge 0)$ is used for
    \inanf{$\psi$ can be {\em inferred or deduced }
    within the {\em calculus } $\KF_S$ augmented by
    {\em premises or additional axioms } $\tup{\varphi}m$ if $m > 0$}.
    We shall write only $\dedq$ for $\ded{S}$ if reference to
    $S$ is clear.
 \end{thms}

 \begin{thms}{Laws of equality}\em\label{eq.laws}
    (Let $a,b,c \in \Expr_{\gamma}$)
    \\
    \begin{tabular*}{.95\textwidth}[t]{cL@{\extracolsep{\fill}}r}
       (I3) &\dedqol{a=b \limp b=a} &(symmetry of equality)
       \\
       (I4) &\dedqol{a=b \limp b=c \limp a=c}
           &(transitivity of equality)
    \end{tabular*}
 \end{thms}
 \begin{proof}
    \begin{mArray}[t]{R@{\qquad}rcc@{\qquad}L}
    (I3) $\Leftarrow$ &%\quad
       a\eqs{\gamma}b
       &\limp (a\eqs{\gamma}a) &\ljnkon{\eqs{\prop}} (b\eqs{\gamma}a)
       &is an inst. of (I2),
    \\
       &&\limp a\eqs{\gamma}a &\limp b\eqs{\gamma}a
       &as $\obl{\eqs{\prop}}$ coincides with $\obl{\leftrightarrow}$
    \\
    &\multicolumn{4}{p{.80\textwidth}}{%
       interchanging the premises by virtue of a
       {\em propositional tautology } and detaching $a=a$ yields (I3)
       }
    \\
    (I4) $\Leftarrow$ &%\quad
       a \eqs{\gamma} b
       &\limp (a \eqs{\gamma} c) &\eqs{\prop} (b \eqs{\gamma} c)
       &is an instance of (I2)
       \\
       &&\limp b = c &\limp a = c
       &as $\obl{\eqs{\prop}}$ coincides with $\obl{\leftrightarrow}$
    \end{mArray}
    \\[-1\baselineskip] \hspace*{.92\textwidth}
 \end{proof}

 \begin{thms}{Equality Theorem}\label{eq.thm}\em
    Let $\obl e \in \Expr \quad
    \obl r, \obl s \in \Expr_{\tau} \quad \vec{y} \in \VAR^*$
    and
    $\obl z \in \VAR_{\tau}$.
    \\ (I5)\Hfill
    If\quad $\gv\,\obl{e} \cap \fv\,\obl{r=s} \subseteq \vcmps{y}$\quad
    then\quad
    $\dedqol{\AQ{\vec{y}}(r=s) \limp e\psub zr = e\psub zs}$
 \end{thms}
 \begin{proof}
       Assume (1) the premises of the theorem,
       (2) \GENP, (3) induction premise, then we observe
       $\obl e = \obl{\opex va}$ and infer the succeeding
    \begin{Enum}[\setcounter{enumi}{3}]
       \item %% (4)
       \begin{math}
          \gv\obl{e}=
          \bigcup\limits_{\indto im}
          (\gv\obl{a_i} \cup \vcmps{v_i})
       \end{math}
       \hfill by \Vref{fbV}{def};
       \item[]
          therefore \quad
          $\gv\obl{a_i} \subseteq \gv\obl{e}$,
          \hfill
          from (1) we obtain
       \item %(5)
          $\gv\obl{a_i} \cap \fv\obl{r=s} \subseteq \vcmps{y}$;
       \item %(6)
          \obl{e\psub zt} =
        \obl{\opexx{v}{a}{\psubvv z{v_i}t{v_i}}}
       \hfill
        by \Vref{subst}{def. }
        \\
       (with $\obl t$ intended to be replaced both by $\obl r$ and $\obl s$)
    \end{Enum}
       \begin{itemize}
       \item[case] $\obl z \in \vcmps{v_i}$:
       \begin{mArray}[t]{rl@{\hspace{4em}}R}
          \obl{a_i[z,\vec{v_i} \gets t,\vec{v_i}]}
          &\multicolumn{2}{l}{= \obl{a_i}
             \hspace*{2.6em} $by \Vref{HS3a}{lemma}$}
          \\
          \obl{a_i[z,\vec{v_i} \gets r,\vec{v_i}]}
          &=
          \obl{a_i[z,\vec{v_i} \gets s,\vec{v_i}]}
          &hence
          \\
          \dedq \mathlq{a_i[z,\vec{v_i} \gets r,\vec{v_i}]}
             &=
             {a_i[z,\vec{v_i} \gets s,\vec{v_i}]}\mathrq
          &(instance of (I1))
       \end{mArray}
       \item[case] $\obl z \notin \vcmps{v_i}$:
       \begin{math}
          \obl{a_i[z,\vec{v_i} \gets t,\vec{v_i}]}
          = \obl{a_i \psub zt}
       \end{math}
       \enspace from \Vref{HS3a}{lemma}
       \\
       \begin{mArray}[t]{lR}
          \dedqol{\AQ{\vec{y}}(r=s) \limp a_i\psub zr = a_i\psub zs}
          &by (3)
          \\
          \vcmps{v_i} \cap \fv\obl{r=s} \subseteq \vcmps{y},
          \quad
          \vcmps{v_i} \cap \fv\obl{\AQ{\vec{y}}(r=s)} = \emptyset
          &from (1)+(4)
          \\
          \dedqol{\AQ{\vec{y}}(r=s) \limp
             \AQ{\vec{v_i}}(a_i\psub zr = a_i\psub zs)}
          &from the preceding 2 lines
       \end{mArray}
       \end{itemize}
       Both cases yield
       \begin{math}
          \dedqol{\AQ{\vec{y}}(r=s) \limp
             \AQ{\vec{v_i}}(%
                a_i\psub{z,\vec{v_i}}{r,\vec{v_i}}
                =
                a_i\psub{z,\vec{v_i}}{s,\vec{v_i}}
                )}
          \vphantom{\Big\vert}
       \end{math}.
       \\[3pt]
       \begingroup %>>>>>>>>>>>>>>>>
       \footnotesize \let\qquad\quad
       Let %\qquad
       \begin{math}%{marray}[t]{l}
          \obl{P_i} = \obl{\Condx vp} \qquad
          \obl{p_i} = \obl{a_i\psub{z,\vec{v_i}}{r,\vec{v_i}}}
          \quad %\\
          \obl{Q_i} = \obl{\Condx vq} \qquad
          \obl{q_i} = \obl{a_i\psub{z,\vec{v_i}}{s,\vec{v_i}}}
       \end{math}%{marray}
       \endgroup %%%<<<<<<<<<<<<<<<<
       \\
       Then% %\quad
       \begin{marray}[t]{l}
          \dedqol{\AQ{\vec{y}}(r=s) \limp \AQ{\vec{v_i}(p_i=q_i)}}
          \\
          \dedqol{\AQ{\vec{v_i}(p_i=q_i)} \limp
             \op(P_1,..,P_{i-1},P_i,Q_{i+1},...)
             =
             \op(P_1,..,P_{i-1},Q_i,Q_{i+1},...)
          }
       \end{marray}
       \\
       \rule{\textwidth}{0pt}
       By chaining these implications and using (I4)
       on p. \vpageref{eq.laws}
       we obtain
       \begin{displaymath}
          \dedqol{\AQ{\vec{y}}(r=s) \limp
             \op(P_1,\dots,P_m) = \op(Q_1,\dots,Q_m)
          }
       \end{displaymath}
       note that
       $\obl{\op(\dots,P_i,\dots)}=\obl{e\psub zr}$
       and
       $\obl{\op(\dots,Q_i,\dots)}=\obl{e\psub zs}$,
       \hfill hence \hfill
       \begin{math}
          \dedqol{\AQ{\vec{y}}(r=s) \limp e\psub zr = e\psub zs}
       \end{math}.
       \quad
       This concludes the induction step.
 \end{proof}

 \begin{korrolar}\em\label{eq.rl}
    (Equality Rule)
    \quad $\obl{r=s} \dedq \obl{e\psub zr = e\psub zs}$
 \end{korrolar}
%
%%####################################################################
%%# funl06.tex  |95-03-12|10:29|                                     #
%%####################################################################
 % FUNL06.TEX %
 % ~~~~~~~~~~~%
 \section{Formalized Theories}
 Throughout the paragraph we assume $S \in \fnlsig$.
 \begin{defin} \em
    \underbar{A functional logic theory }
    is an extension of the {\em calculus } $\KF_S$
    by adding a set $\T \subseteq \Expr_S^{\pi}$ to the axioms.
    The resulting {\em system } is symbolized by $\KF_S[\T]$.
    Members of $\T$ are called {\em nonlogical axioms } of the
    system.
    (We shall often refer to $\KF_S[\T]$ when using $\T$)
 \end{defin}

 \begin{thms}{Notation}\em
       $\T \dedol{S}{\phi}$
       iff there is a {\em proof } of $\obl{\phi}$ within the
       {\em formal system } $\KF_S[\T]$.\quad
       (alias $\phi$ is a {\em theorem of } $\T$,
       or can be {\em inferred from } $\T$)
    \quad \footnotesize
    We shall omit $S$ and write $\T \dedq \obl{\varphi}$ instead
    if only one $S$ is considered. In the sequel we shall always
    assume $\fnlsig S$.
 \end{thms}

 \begin{defin}\label{Df.kons} \em
    \begin{math}
       \kons(\eus{T}) \leqv
       \eus{T} \subseteq \F_{\prop}
       \land  \eus{T} \Ndedq \obl{\falsum}
    \end{math}
 \end{defin}

  \begin{thms}{Compactness Theorem}\label{KompThm}
    \noem
    Each {\em theorem } of $\T \subseteq \Expr_{\pi}$
    is inferable from a {\em finite subset } of $\T$, that is
    \begin{math}
       \T \subseteq \Expr_{\pi} \limp \T\dedqol{\phi} \leqv
       \EQ{\Delta}
       \parenth{big}{\Delta \subseteq \EuScript{T}
          \land \text{finite}(\Delta)
          \land \Delta \dedqol{\phi}
       }
    \end{math}
  \end{thms}
  \begin{proof}
    A (formal) proof only comprises a finite number of axioms,
    hence only a finite subset of $\T$.
  \end{proof}

  \begin{defin}\label{vollst}\label{vollst[]}
    \em %\renewcommand{\S}{\mathrm{S}}
    $\T$ is called \textit{complete relative to} $S$,
    if $\T \subseteq \F^{\prop}_S$ and each closed formula of $S$
    is decidable in $\T$, that is
  \end{defin}
    \skipbasel{-1}
    \begin{displaymath}
       \vollst_S(\eus{T}) \leqv
       \eus{T} \subseteq \F_S^{\prop}
       \land
       \AQu{\obl{\phi}\in\F_S^{\prop}[]}
       (\enspace
          \eus{T}\ded{S}\obl{\phi} \lor \eus{T}\ded{S}\obl{\neg\phi}
       \enspace)
    \end{displaymath}
  (Dependency from $S$ is significant, but we shall
    omit $S$ if confusion is impossible)

  \begin{defin}\noem {$\erw{S_1}{S_2}$}\label{erw}
    symbolizes $(\text{for } \fnlsig S_1,S_2)\colon$
    \inanf{$S_2$ is an extension of $S_1$}, i.e.
    all components of $S_1$ and $S_2$ but $\SOP$ and $\sig$ agree,
    $\SOP_{S_1} \subseteq \SOP_{S_2}$ and
    $\sig_{S_1} = \SOP_{S_1} \restr \sig_{S_2}$.
 \end{defin}

 \begin{thms}{Deduction Theorem}
    \label{dedth}
    \em \enspace (a special form)
    \\ \Hfill
    \begin{math}
       \Gamma\subseteq\F_{\prop}
       \land
       \obl{\phi}\in\F_{\prop}[] \land \obl{\psi}\in\F_{\prop}
       \boldsymbol{\limp}
       \parenth{big}{%
          \Gamma\cup\encurs{\obl{\phi}} \dedq\obl{\psi}
          \boldsymbol{\limp}
          \Gamma \dedq \obl{\phi \limp \psi}
          }
    \end{math}
 \end{thms}
 \begin{proof}
    see \cite{shoenf}, p.33
 \end{proof}

 \begin{thms}{Theorem of Lindenbaum}\label{Lindenbaum}
    For a consistent theory a consistent complete simple extension exists.
 \end{thms}
 \vspace{-1\baselineskip}
 \begin{displaymath}
    \T\subseteq\F_S^{\pi} \land \kons(\T) \limp
    \EQ{\T_2}(\T \subseteq \T_2 \land \kons(\T_2) \land \vollst_S(\T_2))
 \end{displaymath}
 % LINDENB.PF % PROOF of Lindenbaums Lemma %
 % ~~~~~~~~~~
 \begin{proof}^{\EquaNull}
    By \ref{Df.Fnl}\ Def. of $\fnlsig$
    $\SOP_S \cup \VAR_S$ can be well-ordered.
    This implicitly applies to the sets of expressions
    $\F_{\gamma}$ and $\F_{\prop}[]$ (= set of closed formulae), too.
    Hence we may suppose an {\em ordinal enumeration }
    $\enangle{\obl{\phi_{\alpha}}}_{\alpha\in\kappa}$
    of $\F_{\prop}[]$, that is
    \begin{equation}
       \F_{\prop}[] = \clabst{\obl{\phi_{\alpha}}}{\alpha < \kappa}
    \end{equation}
    We claim the following {\it Lemma:}
    \begin{multline}\label{lindenb:lemma}
       \AQ{\eus{T}}\AQ{\obl{\psi}}(\;
       \eus{T}\subseteq\F_{\prop} \land \obl{\psi}\in\F_{\prop}[]
       \land \kons(\eus{T})
       \limp
       \\ \limp
       \kons(\eus{T} \cup \encurs{\obl{\psi}})
       \lor
       \kons(\eus{T} \cup \encurs{\obl{\neg\psi}})
       \;)
    \end{multline}
    For if the opposite is assumed, then
    $\lnot \kons(\eus{T} \cup \encurs{\obl{\psi}})$\quad and
    \\
    $\lnot \kons(\eus{T} \cup \encurs{\obl{\neg\psi}})$.
    By \ref{Df.kons} then
    $\eus{T}\cup\encurs{\obl{\psi}} \dedq \obl{\falsum}$
    and
    $\eus{T}\cup\encurs{\obl{\neg\psi}} \dedq \obl{\falsum}$.
    By \ref{KompThm} (compactness theorem) and the premise
    $\kons(\eus{T})$ we conclude that there are
    $\obltup{\chi}{k}, \obltup{\theta}{\ell}$
    so that
    \begin{math}
       \encurs{\obltup{\chi}{k},\obltup{\theta}{\ell}}
       \subseteq \eus{T}
       \komma \\
       \encurs{\obltup{\chi}{k},\obl{\psi}} \dedq \obl{\falsum}
       \quad \text{and} \quad
       \encurs{\obltup{\theta}{\ell},\obl{\neg\psi}} \dedq \obl{\falsum}
    \end{math}.
    With \ref{dedth} (deduction theorem) we obtain
    \begin{math}
       \encurs{\obltup{\chi}{k}}
       \dedq
       \obl{\psi \limp \falsum}
       \text{ and }
       \\
       \encurs{\obltup{\theta}{\ell}}
       \dedq
       \obl{\neg\psi \limp \falsum}
    \end{math};
    thus
    \begin{math}
       \encurs{\obltup{\chi}{k},\obltup{\theta}{\ell}}
       \dedq \obl{\falsum}
    \end{math},
    that is $\lnot \kons(\eus{T})$ in contradiction to the above premise.
    This confirms (\ref{lindenb:lemma}) {\em Lemma}.
    Using (1) we inductively define an ordinal sequence
    $\enangle{\eus{T}_{\alpha}}_{\alpha\in\kappa+1}$:

    \newcommand{\subeqn}[1]{\tag{\arabic{equation}.#1}}
    \newcommand{\Subeqn}[1]{\mbox{({\arabic{equation}.#1})}\quad}
    \newcommand{\KFT}[1]{\eus{T}_{#1}}
    \refstepcounter{equation}

 \noindent
    \begin{mArray}{rl}
       \Subeqn1
       \eus{T}_0 &= \EAx_T
       \qquad \Subeqn2
       \eus{T}_{\alpha+1} =
       \begin{cases}
          \eus{T}_{\alpha} \cup \encurs{\obl{\phi_{\alpha}}}
          &\text{ if }
          \kons(\eus{T}_{\alpha} \cup \encurs{\obl{\phi_{\alpha}}})
          \\
          \eus{T}_{\alpha} \cup \encurs{\obl{\neg\phi_{\alpha}}}
          &\text{ otherwise }
       \end{cases}
       \\
       \Subeqn3
       \eus{T}_{\lambda} &=
       \textstyle
       \bigcup\limits_{\xi\in\lambda} \eus{T}_{\xi}
       \qquad
       \text{for a limit ordinal } \lambda \in (\kappa+1)
    \end{mArray}
    \refstepcounter{equation}

    we observe
    \\[4pt]
    \begin{mArray}^{\renewcommand{\arraystretch}{1.25}}{r@{\quad}l}
       \Subeqn1&\kons(\KFT{0})
       %\subeqn1
       \hspace{4em} \Subeqn2
       \alpha\in\kappa \limp
       \kons(\KFT{\alpha}) \limp \kons(\KFT{\alpha+1})
       %% \subeqn2
       \\
       \Subeqn3%
       &\text{limit ordinal } \lambda \in (\kappa+1)
       \land \AQu{\alpha\in\lambda} \kons(\KFT{\alpha})
       \limp \kons(%
       \rlap{$\bigcup\limits_{\alpha\in\lambda} \KFT{\alpha})$}
    \end{mArray}
    \\
    {\em Subproof. } If
    $\boldsymbol{\lnot}
       \kons(\bigcup\limits_{\alpha\in\lambda} \KFT{\alpha})
    $,
    then by \ref{KompThm} (compactness theorem) there
    are $\obltup{\theta}{k}$ so that
    \hfill
    \begin{math} %% \begin{equation} \tag{*}
       (*)\enspace
       \obltup{\theta}{k} \in
       \textstyle
       \bigcup\limits_{\alpha\in\lambda} \KFT{\alpha}
       \boldsymbol{\land}
       \encurs{\obltup{\theta}{k}} \dedq \obl{\falsum}
    \end{math}.
    Thus, there are $\tup{\alpha}k < \lambda$, so that
    $\AQu{\indto ik}\; \theta_i \in \eus{T}_{\alpha_i}$.
    Let $\beta=\bigcup\limits_{\indto ik}\alpha_i$ then also
    $\beta < \lambda$ and
    $\encurs{\obltup{\theta}{k}} \subseteq \eus{T}_{\beta}$,
    due to (*) then $\boldsymbol{\lnot} \kons(\eus{T}_{\beta})$
    contradicting the antecedent of implication (4.3). This confirms (4.3).
    By {\em ordinal induction } (limited to range $\kappa+1$)
    we obtain (5);
    By (1) using \ref{vollst[]}{ \em proposition }
    we infer (6) (Subproof below).
    \\[3pt] \nexteq
       $\kons(\eus{T}_{\kappa})$
       \hfill
       \nexteqwl{vollst T_kappa}
       $\vollst(\eus{T}_{\kappa})$
       \hspace{2in}
       \hfill
    \\[3pt]
    {\em Subproof. }
    According to (1) for an arbitrary formula
    $\obl{\psi} \in \F_{\prop}[]$
    $\alpha < \kappa$ exists
    so that $\obl{\psi}=\obl{\phi_{\alpha}}$.
    (3.2) yields
    $\obl{\psi} \in \eus{T}_{\alpha+1} \lor
     \obl{\neg\psi} \in \eus{T}_{\alpha+1}
    $.
    From $\alpha+1 \le \kappa$ by (3.1-3.3) obviously
    $\eus{T}_{\alpha+1} \subseteq \eus{T}_{\kappa}$,
    hence
    $\obl{\psi} \in \eus{T}_{\kappa} \lor
     \obl{\neg\psi} \in \eus{T}_{\kappa}
    $.
    $\obl{\psi} \in \eus{T}_{\kappa}$
    implies $\eus{T}_{\kappa} \dedq \obl{\psi}$, hence
    \begin{math}
       \eus{T}_{\kappa} \dedq \obl{\psi}
       \lor
       \eus{T}_{\kappa} \dedq \obl{\neg\psi}
    \end{math}.
    This confirms (\ref{vollst T_kappa})
 \end{proof}

 \EquaNull
 \begin{theorem}\em\label{Henkin} \textbf{(Henkin)}
    For a theory $\T$ a \textit{conservative extension }
    $\T'$ and a mapping
    \begin{displaymath}%(1)
       \obl{\varphi} \mapsto \obl{c_{\varphi}}
       \colon
       \textstyle\bigcup\limits_{\obl{x} \in \VAR}
       \Expr^{\prop}_{S'}[\obl{x}] \to \COP_{S'}
       \quad \text{exists}
       \qquad (\text{note that}\enspace \COP_{\xi} \subseteq \Expr_{\xi}[]),
    \end{displaymath}
    so that for each
    $\obl{x} \in \VAR$ and $\obl{\varphi} \in \Expr^{\prop}_{S'}[\obl{x}]$:
    \enspace
    \begin{math}
       \obl{(\EQ{x}\,\varphi \limp \varphi\psub{x}{c_{\varphi}})} \in \T'
    \end{math} and
    \\
    $\sig\obl{c_{\varphi}}=\sig\obl{x}$.\quad
    \textit{Conservative extension } means
    \begin{displaymath}%(3)
       \erw{S}{S'} \quad \T \subseteq \Expr^{\pi}_S
       \quad \T' \subseteq \Expr^{\pi}_{S'} \quad \T \subseteq \T'
       \quad\text{and}\quad
       \AQu{\obl{\psi}\in \Expr^{\pi}_{S}}
       \parenth{big}{\T' \dedol{S'}{\psi} \mathbf{\limp} \T \dedol{S}{\psi}}
    \end{displaymath}
 \end{theorem}
 \begin{proof} (as in \cite{shoenf}, p.46)
    Starting with $S_0=S$ we inductively define a sequence of extensions:
    from $S_k$ we obtain $S_{k+1}$ by adding new constant symbols
    $\obl{c_{\varphi}}$, each for one
    \begin{math}
       \obl{\varphi} \in
       \Expr^{\prop}_{S_k}[\obl x] \setminus \Expr^{\prop}_{S_{k-1}},
       \obl{x} \in \VAR
    \end{math}
    (if $k=0$ suppress $\obl{\setminus \Expr^{\prop}_{S_{k-1}}}$)
    and
    \begin{math}
       \sig_{S_{k+1}}\obl{c_{\varphi}} = \sig_S\obl{x}
    \end{math}
    if
    \begin{math}
       \obl{\varphi} \in \Expr^{\prop}_{S_k}[\obl x]
    \end{math}.
    Let $S'$ be the extension of $S$ by adding
    \begin{math}
       \clabst{\obl{c_\varphi}}{%
       \EQu{k}\parenth{big}{%
       \obl{\varphi} \in \Expr^{\prop}_{S_k}[\obl x]
       \land \obl{x} \in \VAR
       }}
    \end{math}
    to component $\SOP$ and stipulating
    \\
    $\sig_{S'}\obl{c_{\varphi}}=\sig_{S}\obl{x}$
    for $\obl{\varphi} \in \Expr^{\prop}_{S_k}[\obl x]$
    and $\obl{x} \in \VAR$;%\quad
    \\ \Hfill
    let
    \begin{math}
       \T' = \T \cup
       \clabst{\obl{(\EQ{x}\,\varphi \limp \varphi\psub{x}{c_{\varphi}})}}%
       {   \obl{\varphi} \in \Expr^{\prop}_{S'}[\obl x]
          \land \obl{x} \in \VAR
       }
    \end{math}.
    \\
    \textit{Claim. } $\T'$ is a \textit{conservative extension}
    of $\T$:
    Assume (1) $\T'\dedol{S'}\varphi$ and
    (2) $\obl{\varphi} \in \Expr^{\prop}_S$.
    Using \Vref{KompThm}{compactness thm.}
    \enspace and\enspace \Vref{dedth}{deduction thm.},
    as $\T' \setminus \T \subseteq \Expr^{\prop}_{S'}[]$,
    we conclude that there are
    $\obl{\psi_1},\dots\obl{\psi_m} \in \T' \setminus \T$
    so that
    \setcounter{equation}{2}
    \begin{equation}^{\jnmusk=3mu}%\tag{3}
       \T \dedqol{\psi_1 \limp \dots \limp \psi_m \limp \varphi}
       \quad \text{and} \quad
       \obl{\psi_1} = \obl{\EQ{x}\theta \limp \theta\psub{x}{c_{\theta}}}
    \end{equation}
    By appropriate arrangement we may assume that $c_{\theta}$
    does not appear in any other $\psi_j\enspace (j=2,..,m)$,
    hence we may replace $c_{\theta}$ by a new variable $y$ so that
    \begin{equation}
       \T \dedqol{(\EQ{x} \theta \limp \theta^x_y)
          \limp \psi_2 \limp \dots \limp \psi_m \limp \varphi
       }
    \end{equation}
    applying well known logical rules and theorems we infer
    \begin{equation}\notag
       \T \dedqol{\EQ{y} (\EQ{x} \theta \limp \theta^x_y)
          \limp \psi_2 \limp \dots \limp \psi_m \limp \varphi}
    \end{equation}
    The antecedent $\EQ{y} (\EQ{x} \theta \limp \theta^x_y)$
    is a logical theorem (deducible from the \textit{variant theorem}
    $\obl{\EQ{x} \theta \limp \EQ{y} \theta^x_y}$),
    hence we can detach it and obtain
    \\
    \begin{math}
       \T \dedqol{\psi_2 \limp \dots \limp \psi_m \limp \varphi}
    \end{math}.
    Iterated application finally yields $\T \dedqol \varphi$ % as required.
 \end{proof}

 \begin{defin} \em \label{sat}
    \enspace For $\FnlSMG{S}{M} \kom \phi \in \Expr_{\prop}$
    and $\eus{T} \subseteq \Expr_{\prop}$
    we define:
 \end{defin}
 \vspace{-1ex}
    \EquaNull
    \begin{eqnarray}
       \M \sat{S} \obl{\phi} &\leqv&
       \AQu{\vec{\sigma} \in \VSRT\stern}
       \AQu{\vec{u} \in \VAR_{\vec{\sigma}}}
       \\
       &&\qquad % \sparenth{big}{%
          \vec{u} \in \persp \obl{\phi}
          \limp
          \AQu{\revec{x} \in \M_{\vec{\sigma}}}
          \M \obl{\phi}(\vec{u})(\revec{x}) = \AlgOp1
       %}
       \notag
       \\
       \M \sat{S} \;\eus{T}\; &\leqv&
       \AQu{\obl{\phi} \in \eus{T}}\; \M \sat{S} \obl{\phi}
       \\
       \eus{T} \sat{S} \obl{\phi} &\leqv&
       \AQ{\M}(\FnlSMG{S}{M} \land \M \sat{S} \eus{T}
          \limp \M \sat{S} \obl{\phi})
    \end{eqnarray}

 \noindent
 For (1) we say
 {$\M$ is a {\em model } of (or $\M$ {\em satisfies})
 $\obl{\varphi}$}, this also applies to (2) w.r.t. $\T$.
 We write $\satq$ for $\sat S$ as long as only one $S$ is considered.
 \begingroup %>>>>>>>>>>>>>>>>
 \renewcommand{\N}{\frak{N}}
 \begin{observ}\em
    If the premises of \fullref{eq.eval} apply to $\M,\N$
    then $\M\sat{S}$ and $\N\sat{S}$ are interchangeable.
    This applies to $\N=\BM$ of \ref{FullStr} as well.
 \end{observ}
 Semantics can be reduced to a notion of {\em full structure }
 as $\BM$, characterized by
 $\BM_{\gamma}^{\vec{\sigma}}=\BM_{\gamma}^{\BM_{\vec{\sigma}}}$
 in place of conditions \FullRef{Strukt}{(3)}.
 \endgroup %%%<<<<<<<<<<<<<<<<

%% file: funl7-bb.tex
% Logic Eprints
%Submitted 0644 Tue Mar 14, 1995 by: a8121dab@helios.edvz.univie.ac.at (josef schoenbrunner )
%logic/schoenbrunner/fun-logic-completeness
%logic/schoenbrunner/fun-logic-completeness/funl7-bb.tex
%

% FUNL7-BB = funl07.tex + funlk.bbl $1
%%####################################################################
%%# funl07.tex  º95-03-12º21:04º                                     #
%%####################################################################
 % FUNL07.TEX %
 % ~~~~~~~~~~~%
 \begingroup %>>>>>>>>>>>>>>>>
 %% FUNL07.DEF %%%%%%% Macro Defs. fr funl07.tex %
 %% ~~~~~~~~~~~~~~~~~~
    \newif\ifmiteqlemma \miteqlemmatrue
    \let\phi\varphi
    \let\norm\underline
    \let\nor\underline
 \newcommand{\inCM}{_{\CM}}
 \newlength{\mathextraskip} \setlength{\mathextraskip}{1ex}
 \allowdisplaybreaks
 \newenvironment{fitbox}[1]{\hbox to .#1\textwidth\bgroup}{\hss\egroup}
 \newcommand{\NL}{\\ \hbox{}&}
 \newcommand{\NLNT}{\\ \notag\hbox{}&}
 \newcommand{\enumlit}[2]{#1^{\text{#2}}}
 \newcommand{\subpf}{{\em Subproof. }}
 \newcommand{\fos}{\Expr_{S}^{\prop}}
 \newcommand{\foss}{\Expr_{S'}^{\prop}}
 %------------------------------------------------%
 \section{Obtaining a Model of a Consistent Theory}

 \begin{thms}{Extension Theorem}\label{ErwThm}\em
    Assume $\T$ is a consistent theory of signature $S$
    (technically
    $\fnlsigs \land \T \subseteq \F^{\prop}_S \land \kons(\T)$),
    then  a so called {\em complete and consistent Henkin theory }
    $\T'$ exists, which is an {extension of } $\T$.
    That is, in the sequel we shall rely on the following conditions:
 \end{thms}
    (1)
       \begin{math}
          \T \subseteq \T'
          \komma\qquad
          \T \subseteq \F^{\prop}_S \subset \F^{\prop}_{S'}
          \supseteq \T'
       \end{math}
       \quad
    (2)\enspace $\erw{S}{S'}$,
       the components $\SRT_{S'}, \VSRT_{S'}, \ISRT_{S'}$
       of signature $S'$ agree with those of $S$, only $\SOP$ and
       $\COP$ differ.
       (In the sequel we shall omit the subscripts for those
        components that agree)
       \quad
    (3)\enspace $\vollst_{S'}(\T')$.
       \\[3pt]
    (4)   There is a map
       \mbox{%-------------------------------------%
       $%\begin{math}
          \ophi \mapsto \obl{\sco_{\phi}} \enspace \boldsymbol{:} \enspace
          \textstyle
          \bigcup\limits_{\obl{x}\in\VAR} \F^{\prop}_{S'}[\obl{x}]
          \; \rightarrow \;       \COP_{S'}
       $%\end{math}
       \enspace}%----------------------------------%
       that assigns to each $\ophi \in \F^{\prop}_{S'}[\obl{u}]$
       (formula with at least one free variable $\obl{u}$)
       the so called {\em special constant } $\obl{\sco_\phi}$
       so that
       \mbox{%-------------------------------------%
       $\T'\dedq \obl{%
             \EQ{^\alfa u} \phi \limp \sub{\phi}u{\sco_{\varphi}}
             }
       $\enspace}%----------------------------------%
       and $\usig_{S'} \obl{u} =
       \usig_{S'} \obl{\sco_{\varphi}} = \obl{\alpha}
       $
 \\[-1ex]
 \begin{pf}
    Use \vref{Henkin}, then for the resulting $\T'$ apply \ref{Lindenbaum}.
    The second extension does not cancel the qualities achieved with the first.
 \end{pf}

 \ifmiteqlemma %:::::::::::%
    Now we consider three consequences of \ref{ErwThm} that shall be
    prerequisites for investigations related to a term-structure
    built upon $S'$ and $\T'$.
 \fi %::::::::::::::::%

 \begin{lemma}\label{HT1}\noem
    If $\T'$ is a {\em Henkin theory } due to \ref{ErwThm}(4),
    \begin{math}
       \varphi\in\Expr_{\pi}[\vec{z}] \komma
       \vec{z} \in \VAR_{\vec{\varrho}} \komma
       \vec{\varrho} \in \VSRT^r
    \end{math}
    and all components of $\vec{z}$ are different
    then there are special constants
    $\obl{c_i} \in \COP_{\varrho_i} % \enspace
     \\ (\indto ir)$
    so that for $\vec{c}=\enangle{c_i}_{\indto ir}$:\qquad
    $\T'\dedqol{\varphi\psubv zc \limp \AQvecz\varphi}$
 \end{lemma}

 \begin{pf}^{\jnmusk=0mu}
    For $\obl{\varphi} \in \Expr_{\prop}[\obl{x}]\komma
       \eqobl{d}{c_{\lnot \varphi}}
    $
    by \ref{ErwThm}(4) % we obtain
    $\T' \dedqol{\varphi\psub xd \limp  \AQ{x}\varphi}$.
    Instances of such theorems are
    $\obl{(\psi_0\limp\psi_1)},\dots,\obl{(\psi_{r-1}\limp\psi_r)}$,
    where each $\eqobl{\psi_i}%
       {\AQ{z_1}...\AQ{z_i}%
       (\varphi\psub{z_{i+1}}{c_{i+1}}...\psub{z_r}{c_r})}$.
    Transitivity of {\em implication } yields
    \\
    $\T'\dedqol{\psi_0\limp\psi_r}$, according to {\em lemma } \fullref{HS5}
    $\eqobl{\psi_0}{\varphi\psub{z_1}{c_1}...\psub{z_r}{c_r}}=
    \\
     \obl{\varphi\psubv zc}$, hence
    $\T'\dedqol{\varphi\psubv zc \limp \AQvecz\varphi}$
 \end{pf}

 \ifmiteqlemma %:::::::::::%
 \begin{lemma}\noem\label{eqlemma1}
    Assume $\T'$\ due to \ref{ErwThm}(4),\spa
    $\vec{\delta}\in\VSRT^{\ell}\spa \vec{u},\vec{v}\in\VAR_{\vec{\delta}}
       \spa\obl{a}\in\Expr_{\gamma}[\vec{u}]
       \\
       \spa\obl{b}\in\Expr_{\gamma}[\vec{v}]
    $ and no component of $\vec{u}$ or $\vec{v}$ occurs more than once.
    \\
    If $\AQu{\vec{d}\in\Expr_{\vec{\delta}}[]}
       \T'\dedqol{a\psubv ud = b\psubv vd}
    $\enspace
    and if $\vec{z}=\enangle{\obl{z_i}}_{\indto{i}{\ell}}
       \in\VAR_{\vec{\delta}}$
    is such that each $\obl{z_i}$ is different from the members of
    $\vec{u}$ and $\vec{v}$ and neither occurs in $\obl a$
    nor in $\obl b$, then \enspace
    $\T'\dedqol{\AQvecz a\psubv uz = b\psubv vz}$
 \end{lemma}
 \begin{pf}% Lemma \ref{eqlemma1}
    Let $\obl{s}=\obl{a\psubv{u}{z}}$
    and $\obl{t}=\obl{b\psubv{v}{z}}$.
    Our supposition for $\vec{z}$ implies
    that $\obl{\psubv uz \psubv zd}$ operates exactly like
    $\obl{\psubv ud}$, hence
    $\obl{s\psubv{z}{d}}=\obl{a\psubv{u}{d}}$ and
    $\obl{t\psubv{z}{d}}=\obl{b\psubv{v}{d}}$.
    Premise $\T'\dedqol{a\psubv ud = b\psubv vd}$ then yields
    $\T'\dedqol{t\psubv{z}{d}=s\psubv{z}{d}}$,
    by \ref{subst} {\em def. } this is (*)
    $%
     \T'\dedqol{(t=s)\psubv{z}{d}}$
    for arbitrary $\vec{d}\in\Expr_{\vec{\delta}}$.
    Now we are ready to use \ref{HT1}, the {\em lemma } that provides
    {\em special constants } $\vec{c}=\enangle{\obl{c_j}}_{\indto{j}{\ell}}$
    (each $\obl{c_j}\in\COP_{\delta_j}$) so that
    $\T'\dedqol{(t=s)\psubv{z}{c} \limp \AQvecz\; t=s}$,
    detachment with (*) yields
    $\T'\dedqol{\AQvecz\; t=s}$,
    recalling the definition of $\obl{s}$ and $\obl{t}$ this is
    $\T'\dedqol{\AQvecz\ b\psubv{v}{z} = a\psubv{u}{z}}$
 \end{pf}

 \begin{lemma}\noem\label{eqlemma2}
    Assume $\T'$ due to \ref{ErwThm}(4), $\pGP$ and $\avar{\pGP}$
    with $\vec{u}=\vec{\avar{u}}=\emptyseq$.
    $\avar{\pGP}$ denotes the change of all symbolic parameters
    except $\oblt{op}, m$ and $r_i$ by appending a star-superscript
    performed on $\pGP$ (as new names for corresponding but
    possibly different things are required).
    \\
    If $\AQu{\indto im}\AQu{\vec{s}\in\Expr_{\vec{\beta_i}}[]}
       \T'\dedqol{\avar{a}_i\psubv{\avar{v}_i}s = a_i\psubv{v_i}s}
    $
    then $\T'\dedqol{\avar{e} = e}$
 \end{lemma}
 \begin{pf}
    Applying the preceding {\em lemma\/} to the last premise yields
    \quad
    \begin{math}
       \AQu{\indto im\hspace{-.3em}}
       \\
       \T'\dedqol{\AQvecz
          \avar{a_i}\psubv{\avar{v_i}}{z_i} =
          a_i\psubv{v_i}{z_i}
          }
    \end{math}
    (by appropriate choice of $\vec{z_i}\in \VAR_{\vec{\beta_i}}$).
    Now we are ready to use {\em equality axioms }
    \ifallowpageref (page \pageref{eqax})\fi.
    Iterated application of (I2) yields the chain of equations
    \begin{math}^{\small}
       \op(P_1,P_2,\dots,P_m)=\op(Q_1,P_2,\dots,P_m)
       =\dots =\op(Q_1,\dots,Q_m)
    \end{math}
    where {\small$\obl{P_i}=\obl{\Condx{\avar{v}}{\avar{a}}}$}
    and {\small$\obl{Q_i}=\obl{\Condx va}$}.
    From (I4) \textit{transitivity} we obtain $\T'\dedqol{L=R}$ with
    {\small
    $\obl{L}=\obl{\op(\tup Pm)}=\obl{\opex{\avar{v}}{\avar{a}}}
    =\obl{\avar{e}}
    $
    and
    $\obl{R}=\obl{\op(\tup Qm)}=\obl{\opex va}
       \undertext{\ref{subst}}=
       \obl{\opex va}=\obl{e}
    $},
    hence $\T'\dedqol{\avar{e}=e}$.
 \end{pf}
 \fi %::::::::::::::::%

 \begin{defin}\em\label{norm} {(norm of closed expressions)}
    \quad Assume $\fnlsig S' \komma \T' \subseteq \F_{S'}^{\pi}$.
    \\
    For each $\gamma \in \SRT$ and
    $\obl e \in \F^{\gamma}_{S'}[] \aliasgl \Expr_{\gamma}[]$
    the {\em norm } $\obl{\norm e}$ is defined as follows:
 \end{defin}
    If  $\gamma = \prop$  then if $\T'\dedq \obl{e}$
    then  $\obl{\norm e}=\obl{\verum}$
    otherwise  $\obl{\norm e}=\obl{\falsum}$.
    \\
    If $\gamma\neq\prop$ then
          we rely on an ordinal enumeration of
          $\Expr_{\gamma}[]$ and define
          $\obl{\norm{e}}$ to be the first within the subset
          \begin{math}
             \clabst{\obl{a}\in\Expr_{\gamma}[]}%
             {\T'\dedq \obl{a \eqs{\gamma} e}}.
          \end{math}
          w.r.t. this enumeration.

 \begin{eigen} \em \label{norm.p}
    Let $\fnlsig S' \komma \T'\subseteq \foss \komma
    \gamma \in \SRT \setminus \encurs{\prop}$
    and $\obl{e},\obl{e_1},\wopts0{\\}\obl{e_2} \in \Expr_{\gamma}[]$
    then
    \begin{math}
    \begin{Array}^{%
          \newcommand{\Hexpr}{\Expr_{\gamma}[]}
          \newcommand{\DED}{\T'\dedq}
          \renewcommand{\arraystretch}{1.5}
       }[t]{l}
       (1) \quad \DED \obl{\norm{e} \eqs{\gamma} e}
       \hspace{3em}
       (2) \quad \DED \obl{e_1 \eqs{\gamma} e_2}
          \leqv \obl{\norm{e_1}} = \obl{\norm{e_2}}
       \\
       (3) \quad \obl{\varphi} \in \Expr_{\prop}[]
          \land \vollst_{S'}(\T')
          \limp
          \DED \obl{\varphi} \leqv
             \obl{\norm{\varphi}} = \obl{\verum}
    \end{Array}
    \end{math}
 \end{eigen}
 \begin{proof} immediately from the preceding definition.
 \end{proof}

    \renewcommand{\a}{\frak{a}}
    \newcommand{\va}{\vec{\a}}

 \begin{defin}\noem\label{TStr}
    {\bf(term structure)}
    Let $\HT$ be a {\em complete and consistent Henkin theory},
    i.e. we presuppose the conditions of
    \ref{ErwThm} except (1) for $\HT$.
    We define the {\em term structure} $\CM$ as a function on
    $\SRT \cup (\SRT \x \VSRT^*) \cup \SOP$.
    {\small
    As we now consider only one {\em signature }
    ($S'$, that of $\HT$ and do not refer to $\T$ and $S$)
    we shall omit subscript $\HS$.
    }
 \end{defin}
 \vspace{-1\baselineskip}
    \EquaNull  %% \let\bsparenth\sparenth
 \begin{align}^{\renewcommand{\arraystretch}{1.4}}
    \hbox{}&\sparenth{normal}{\gamma \in \SRT} \qquad %(1)
    \CM(\gamma) \aliasgl \CM_{\gamma} =
    \clabst{\obl{\nor e}}{\obl e \in \Expr_{\gamma}[]}
    \\%[3pt]
    \hbox{}&%
    \begin{Array}[t]{rl} %(2)
       \hbox{}&\sparenth{big}{\gamma \in \SRT \komma
       \vec{\sigma} \in \VSRT^*
       } \quad \CM(\gamma,\vec{\sigma})
       \aliasgl
       \CM_{\gamma}^{\vec{\sigma}} \boldsymbol{=}
       \\
       &\boldsymbol{=}
       \begin{cases}
          \CM_{\gamma}^{\enangle{}} = \CM_{\gamma}
          &\text{if } l = 0
          \\
          \boldsymbol{\lbrace} \msf{f} \boldsymbol{\mid}
             \EQu{\vec{u} \in \VAR_{\vec{\sigma}}}
             \EQu{\obl{e} \in \Expr_{\gamma}[\vec{u}]}
             \msf{f} =
             \funkd{\CM(\vec{\sigma})}{\CM(\gamma)}%
             {\vec{c}}{\obl{\nor{\sub{e}{\vec{u}}{\vec{c}}}}}
          \boldsymbol{\rbrace}
          &\text{if } l\neq0
       \end{cases}
    \end{Array}
    \\[3pt]\hbox{}&%
       \sparenth{big}{\oblt{op} \in \SOP \komma
          \sigoop = \stdsig
       }
       \quad
    \hbox to .2\textwidth{$%
          \text{If } m = 0 \quad
          \CM\oblt{op} = \oblt{\underbar{op}}
          $\hss}
       \NLNT
       \text{If } m \neq 0 \quad \CM\oblt{op} =
    \hbox to .5\textwidth{%
       $\clabst{\opaar{\va}{\obl{\nor{e}}}}{\enspace
          \va = \enangle{\a_i}_{\indto im} \in
             \Prod_{\indto im} \CM_{\alpha_i}^{\vec{\beta_i}}
          \land
          \EQ{\vec{a}}\EQ{\VVec{v}} \text{(i)-(iv)}
       }$%
          \hss}
       \NLNT
    \hbox to .6\textwidth{%
          (i)\quad
          $\vec{a} = \enangle{\obl{a_i}}_i %{\indto im}
          \in \Prod_{\indto im} \Expr_{\alpha_i}[\vec{v_i}]
          $\qquad (ii)\quad
          $\VVec{v} = \enangle{\vec{v_i}}_i %{\indto im}
          \in \Prod_{\indto im} \VAR_{\vec{\beta_i}}
          $
          \hss}
       \NLNT
    \hbox to .6\textwidth{%
          (iii)\quad
          $\AQu{\indto im} \enspace \frak{a}_i =
          \begin{cases}
             \obl{\nor{a_i}} &\text{ if } r_i = 0
             \\
             \funkd{\CM(\vec{\beta_i})}{\CM(\alpha_i)}%
             {\vec{b}}{\obl{\nor{\sub{a_i}{\vec{v_i}}{\vec{b}}}}}
             &\text{ if } r_i \neq 0
          \end{cases}
          $
          \hss}
       \\[1ex]\hbox{}&\notag
       \text{(iv)}\quad\obl{e} = \obl{\opex va}
 \end{align}

 \begin{remark}\em\label{imp.pGP}
    (i)-(iv) in \ref{TStr}(3) imply $\pGP$ with $\vec{u}=\emptyseq$
 \end{remark}

 \begin{eigen}\em\label{TM:Strukt}
    \qquad $\FnlSMG SX$
    \qquad (Def. s. \fullref{Strukt})
 \end{eigen}
 \noindent
 % PROOF71.TEX
 % ~~~~~~~~~~~
 \begin{proof}
    We show that the laws of Def. \fullref{Strukt} apply to $\CM$
    (in place of $\M$)
    referring to them with the numeration used in that definition.
    We rely on {\em lemmas } \ref{HS1} to \ref{HS5} and
    Def. \fullref{TStr}\ of $\CM$.
    For (1) and (3.1) there is nothing to prove.
    Our first goal is to prove validity of (2), this only requires
    verification that $\CM\,\oblt{op}$ is a function if $m\neq 0$.
    We refer to \ref{TStr}(3), case $m\neq 0$ (def. of $\CM\,\oblt{op}$).
    The {\em goal } is then reduced to the task of deduction
    from the {\em premises }
    $\opaar{\vec{\frak{a}}}{\obl{\nor e}} \in \CM\,\oblt{op}
       \komma
       \opaar{\vec{\avar{\frak{a}}}}{\obl{\nor{\avar{e}}}} \in \CM\,\oblt{op}
    $
    and $\vec{\avar{\frak{a}}} = \vec{\frak{a}}$
    to the {\em conclusion } $\obl{\nor e}=\obl{\nor{\avar e}}$.
 {\em Subproof (Sketch).}
    By expansion of $\enumlit1{st}$ and $\enumlit2{nd}$ premise
    according to \ref{TStr}(3)(case $m\neq 0$) we obtain
    (i)-(iv) of \ref{TStr}(3) and a starred isomorphic variant
    $\avartext{(i)}-\avartext{(iv)}$. By \ref{imp.pGP} {\em remark }
    we succeed to $\pSyntass$ and $\avar{\pSyntass}$ with
    $\vec{u}=\vec{\avar u}=\emptyseq$.
    Still one condition lacks
    to use \ref{eqlemma2} {\em lemma}.
    From {\em premise } $\vec{\avar{\frak{a}}} = \vec{\frak{a}}$
    and expanding (iii) of \ref{TStr}(3) we obtain
    {\em equations }
    $%
       \obl{\nor{\avar{a_i}\psubv{\avar{v_i}}{b}}} =
       \obl{\nor{a_i\psubv{v_i}{b}}}
    $
    for each $\vec{b}\in\CM_{\vec{\beta_i}}$.
    These {\em equations } can be transformed into theorems of $\T'$
    according to \ref{norm.p}(2) and then generalized by \ref{norm.p}(1)
    to apply to each $\vec{b} \in \Expr_{\vec{\beta_i}}[]$,
    so that finally   \ref{eqlemma2} can be used to deduce
    $\T'\dedqol{\avar{e}=e}$. To present this in detail:
    \ref{norm.p}(2) yields
    $\T'\dedqol{\avar{a_i}\psubv{\avar{v_i}}{b}= a_i\psubv{v_i}{b}}$
    for $\vec{b}\in\CM_{\vec{\beta_i}}$.
    Taking into account, that
    $\obl{d_j}\in\Expr_{\beta\biind ij}[]$
    implies $\obl{\nor{d_j}}\in\CM_{\beta\biind ij}$
    we obtain
    $\AQu{\vec{d}\in\Expr_{\vec{\beta_i}}\hspace{-.27em}}
    \T'\dedqol{\avar{a_i}\psubv{\avar{v_i}}{\nor d}= a_i\psubv{v_i}{\nor d}}
    $.
    Referring to \ref{norm.p}(1): $\T'\dedqol{\nor{d_j}=d_j}$
    by virtue of \Vref{eq.rl}{corollary (equality rule)}
    we may replace each $\obl{\nor{d_j}}$ by $\obl{d_j}$,
    hence $\vec{\nor{d}}$ by $\vec{d}$
    and are ready to apply \ref{eqlemma2}, that yields
    $\T'\dedqol{\avar{e}=e}$
    and by \ref{norm.p}(2) $\obl{\nor{\avar{e}}}=\obl{\nor{e}}$ \qed(2).
    We turn to the {\em completion qualities}. Next goal is (3.2 a)
    \begin{math}
       \text{cst}^{\vec{\sigma}}_{\reob{w}} \in \CM_\gamma^{\vec{\sigma}}
       \quad
       \text{for } \reob{w} \in \CM_\gamma
    \end{math}.
    \subpf
    Suppose $\reob{w} \in \CM_\gamma$,
    then (applying (1) of \ref{TStr})
    $\reob w = \obl{\nor e}$ for some $\obl e \in \Expr_\gamma[]$,
    by \fullref{frV+pers} $\fv\obl{e}=\emptyset$, hence by \fullref{HS1}
    $\obl{e} = \obl{\sub{e}{\vec{u}}{\vec{c}}}$,
    then
    $%
       \text{cst}^{\vec{\sigma}}_{\reob{w}} =
       \text{cst}^{\vec{\sigma}}_{\obl{\nor e}} =
       \funkd{\CM_{\vec{\sigma}}}{\CM_{\gamma}}%
       {\vec{c}}{\obl{\nor{\sub{e}{\vec{u}}{\vec{c}}}}}
    $,
    that is
    $\text{cst}^{\vec{\sigma}}_{\reob{w}} \in \CM_{\gamma}^{\vec{\sigma}}$
    due to (2) of \ref{TStr} Def.
    \qed(3.2 a). %% This confirms (3.2 a) \qed.
    Next goal is (3.2 b)
    \begin{math}
          \text{pj}^{\vec{\sigma}}_j
          =
          \funkd{\CM\vec{\sigma}}{\CM\sigma_j}{\vec{c}}{\obl{c_j}}
                \in\CM_{\sigma_j}^{\vec{\sigma}}
    \end{math}.
    \\
    \subpf
    Let $\vec{c}=\enangle{\obl{c_i}}_{\indto il}$.
    If we apply \ref{subst} Def.,
    case $m=0 \kom \oblt{op} \in \vcmps x
    %% \clabst{\obl{x_i}}{\indto il}
    $,
    perform a substitution such that the second
    {\em condition of case } changes into
    $\obl{u_j} \in \vcmps u$ %%%% \clabst{\obl{u_i}}{i\dots}$
    and use $\obl{c_j}=\obl{\nor{c_j}}$
    (as $\obl{c_j} \in \CM_{\sigma_j}
     \subseteq \clabst{\obl{\nor e}}{\obl{e} \in \Expr_{\sigma_j}}$),
    then we obtain
    \begin{math}
       \text{pj}_j^{\vec{\sigma}} =
        \funkd{\CM{\vec{\sigma}}}{\CM{\sigma_j}}%
       {\vec{c}}{\obl{\nor{\sub{u_j}{\vec{u}}{\vec{c}}}}}.
    \end{math}
    This matches the pattern of \ref{TStr} (2)
    for $\text{pj}_j^{\vec{\sigma}} \in \CM_{\sigma_j}^{\vec{\sigma}}$
    \qed(3.2 b). %% This confirms (3.2 b) \qed.
    Next goal is
    \\(3.3)
    \begin{math}
       \vec{c} = \enangle{\obl{c_i}}_i \in \CM_{\vec{\sigma}}
       \land
       \reob{g}\in\CM_\gamma^{\vec{\sigma}\verkett\vec{\varrho}}
       \boldsymbol{\limp} \reob{g}_{\vec{c}}
       \in \CM_{\gamma}^{\vec{\varrho}}.
    \end{math}
    \subpf
    Applying \ref{TStr} \wopts1{ (Def. of $\CM$) }
    we can replace the $2^{\text{nd}}$ premise by
    $\reob{g}=
       \funkd{\CM{\vec{\sigma}\vec{\varrho}}}{\CM{\gamma}}%
       {\vec{c}\vec{d}}{\obl{\sub{e}{\vec{x}\vec{y}}{\vec{c}\vec{d}}}}
    $.
    If we replace $\reob g$ in (3.3) by the right hand term,
    it only remains to show:
    if $\vec{x}\in\VAR_{\vec{\sigma}} \kom
       \vec{y}\in\VAR_{\vec{\varrho}}
       \kom \obl{e}\in \Expr_{\gamma}[\vec{x}\vec{y}]
    $ and
    $\vec{c}=\enangle{\obl{c_i}}_{\indto{i}{\ell(\sigma)}}
     \in \CM_{\vec{\sigma}}
    $, then
    $\reob{g}_{\vec{c}} =
       \funkd{\CM{\vec{\varrho}}}{\CM{\gamma}}%
       {\vec{d}}{\obl{\nor{e\psubvv xycd}}}
       \in \CM_{\gamma}^{\vec{\varrho}}
    $.
    The expression for $\reob{g}_{\vec{c}}$ must be transformed in
    such a manner that it matches the pattern of $\reob f$
    in \FullRef{TStr}{(2)}.
    This is achieved by aid of {\em lemma } \ref{HS3} :
    $\obl{e\psubvv xycd} = \obl{e\psubvv xycy \psubv yd}$.
    Substitution yields
    $\reob{g}_{\vec{c}} =
       \funkd{\CM{\vec{\varrho}}}{\CM{\gamma}}%
       {\vec{d}}{\obl{\nor{e\psubvv xycy \psubv yd}}}
    $.
    According to \ref{TStr}\ (2) then
    $\reob{g}_{\vec{c}}\in\CM_{\gamma}^{\vec{\varrho}}$ \qed (3.3).

    Next we show (3.4).
    {\em Outline: } It requires inheritance of the quality
    $\CM^{\vec{\sigma}}_{\gamma}$ of a map
    $\CM_{\vec{\sigma}}\to\CM_{\gamma}$,
    when composition of such maps with $\CM\oblt{op}$
    due to \ref{Strukt}{ \em definition }(3.4) \apageref{Strukt}
    is performed.
    Quality $\CM^{\vec{\sigma}}_{\gamma}$ applies to
    $\reob h \colon \CM_{\vec{\sigma}}\to\CM_{\gamma}$ iff
    $\reob h$ can be defined by an expression $\obl{e}$
    as a substitution mapping $\vec{c}\mapsto\obl{\nor{e\psubv xc}}$.
    To show this {\em inheritance } we start from \ref{TStr}(3),
    compare the terms that arise with those
    obtained from quality $\CM^{\vec{\sigma}\vec{v_i}}_{\alpha_i}$
    applied to $\reob h_i \vec{c}$.
    The goal that the composed map
    $\vec{c}\mapsto\CM\oblt{op}(\dots,\reob{h}_i\vec{c},\dots)$
    belongs to $\CM^{\vec{\sigma}}_{\gamma}$
    is achieved if we find an $\obl e \in \Expr_{\gamma}$ so that
    $\CM\oblt{op}(\dots,\reob{h}_i\vec{c},\dots)=\obl{\nor{e\psubv xc}}$.
    The solution will be quite natural $\obl{e}=\obl{\opex va}$.

    {\em Subproof of } (3.4).
    From the premises (p1)-(p4) we infer (|) below.
    \\
    \EquaNull
    \begingroup %>>>>>>>>>>>>>>>>
    \renewcommand{\theequation}{p\arabic{equation}}
    \begin{math}
       \nexteq m>0 \komma \ell>0 \hspace{1in} %(p1)
       \nexteqwl{p2}
       \oblt{op} \in \SOP
       \komma
       \sig \oblt{op} = \stdsig
    \end{math}
    \\
    Letter $i$ is understood to be bound by $\AQu{\indto im} \dots$
    in subsequent context.
    \\
    \begin{math}
       \nexteqwl{p3}
          \reob{g}_i \in \CM_{\alpha_i}^{\vec{\sigma}\vec{\beta_i}},
          \enspace \text{that is owing to \ref{TStr} (2), case $\ell>0$ :}
          \\ \hbox{\quad} \hfill %\enspace
          \reob{g}_i =
          \funkd{\CM\vec{\sigma}\vec{\beta_i}}{\CM\alpha_i}%
          {\vec{c}\vec{\beta_i}}%
          {\obl{\nor{a_i\psubvv{x}{v_i}{c}{b_i}}}}
          \text{ with } a_i \in \Expr_{\alpha_i}[\vec{x}\vec{v_i}]
          \\ % \quad
       \nexteqwl{p4}
          \text{ if } r_i=0 \text{ then } \reob{h}_i = \reob{g}_i
             \text{ else } \reob{h}_i =
             \funkd{\CM_{\vec{\sigma}}}{\CM_{\alpha_i}^{\vec{\beta_i}}}%
             {\vec{c}}{{\reob{g}_i}_{\vec{c}}}
          \text{ where } {\reob{g}_i}_{\vec{c}} =
          \funkd{\CM\beta_i}{\CM{\alpha_i}}{\vec{b_i}}%
          {\obl{\nor{a_i\psubvv{x}{v_i}{c}{b_i}}}}.
    \end{math}
    \\
    Our goal is to show (|)
    $%
       \funkd{\CM\vec{\sigma}}{\CM\gamma}{\vec{c}}%
       {\CM\oblt{op}(\dots,\reob{h}_i\vec{c},\dots)}
       \in
       \CM_{\gamma}^{\vec{\sigma}}
    $.
    We now turn to conclusions from (p1-p4).
    \\[1ex]
    \nexteqwl{p5}
    $\obl{a_i\psubvv{x}{v_i}{c}{b_i}} =
       \obl{a_i\psubvv{x}{v_i}{c}{v_i}\psubv{v_i}{b_i}}
    $
    (according to  \ref{HS3}\ {\em lemma}).
    \\[1ex]
       From \ref{p4}+\ref{p5} we infer
    \\
    \begin{math}
    \nexteqwl{p6}
       \reob{h}_i \vec{c} =
       \begin{cases}
          \reob{g}_i \vec{c} = \obl{\nor{a_i\psubv xc}}
          &\text{ if } r_i=0
          \\
          {\reob{g}_i}_{\vec{c}} \colon
          {\vec{b_i} \mapsto
             \obl{\nor{a_i\psubvv{x}{v_i}{c}{b_i}\psubv{v_i}{b_i}}}
             }\quad
          &\text{ if } r_i \neq 0
       \end{cases}
    \end{math}

    \noindent
    In order to verify (|) we look at \ref{TStr} (2) Def. $\CM$.
    It requires transforming of
    $\CM\oblt{op}(..., \reob{h}_i \vec{c}, ...)$
    into an expression of shape $\obl{\nor{\op(...,?,...)}}$.
    This can be achieved by application of \ref{TStr} (3)
    if we find a substitute of `?' so that
    $\CM\oblt{op}(..., \reob{h}_i \vec{c}, ...) = \obl{\nor{\op(...?...)}}$.
    We have to show (||)
        $\opaar{\enangle{\reob{h}_i \vec{c}}_{\indto im}}%
          {\obl{\nor{\op(...?...)}}}
       \in \CM\oblt{op}
    $ for some `?'.
    We observe the following correspondence of terms in \ref{TStr} (3)
    and those from our current {\em Situation}:

    \noindent
    \begin{tabular}{rCCCC}
       &&(\text{if } r_i=0) &(\text{if } r_i>0)
       \\ \hline
       \ref{TStr}(3):\enspace &\vec{\frak{a}}
       &\obl{\nor{a_i}}
       &\obl{\nor{a_i\psubv{v_i}{b}}}
       &\obl{\opex va}
       \\
       {\em Situation:\enspace} &\reob{h}_i \vec{c} &\obl{\nor{a_i\psubv xc}}
       &\obl{\nor{a_i\psubvv x{v_i}c{b_i}}}
       &\obl{\op(\dots,?,\dots)}
       \\[3pt]
       \hline
    \end{tabular}
    \\[1ex]
    This makes evident what to be substituted for `?' in (||): only
    $\smash[b]{\obl{\Condx va \psubv xc}}$ may be considered.
    The problem is now reduced to showing
    \begin{equation}\tag{|||}
       \hspace*{-3em}
       (\text{for } r_i>0\quad
          \vec{c}\in\CM_{\vec{\sigma}} \quad
          \vec{b_i} \in \CM_{\vec{\beta_i}}
       ) \quad
       \obl{\nor{a_i\psubvv{x}{v_i}{c}{b_i}}} =
       \rlap{$
       \obl{\nor{a_i\psubv{x}{c} \psubv{v_i}{b_i}}}
       $}
    \end{equation}
    \noindent
    As $\CM_{\vec{\sigma}} \subseteq \Expr_{\sigma}[]$,
    this equation is supplied by \ref{HS5}\ {\em lemma},
    but only if $\vcmps{x} \cap \vcmps{{v_i}} = \emptyset$.
    We can achieve this requirement by renaming the variables
    $\vec{v_i}$ using {\em equality axioms } (I2) of \fullref{CAL},
    namely the following instance:
    \begin{multline}\label{p7}
       \AQ{\vec{w_i}}
       \;
       a_i\psubv{v_i}{w_i} =
       a_i\psubv{v_i}{q_i}\psubv{q_i}{w_i}
       \limp
       \\
       \limp
       \op(\Delta,\Condx va,\Gamma) =
       \op(\Delta,\Condx qa\psubv{v_i}{q_i},\Gamma)
    \end{multline}
    (in place of $\obl{b_i}$ at % \marginpar{at correct?}
    \ref{CAL}\ now $\obl{a_i\psubv{v_i}{q_i}}$).
     If we choose $\vec{w_i}$ and $\vec{q_i}$ so that
     the variables of the sequences
     $\vec{v_i}\kom\vec{w_i}$ and $\vec{q_i}$ are pairwise disjoint,
     then also\\
     $\obl{a_i\psubv{v_i}{w_i}} =
      \obl{a_i\psubv{v_i}{q_i}\psubv{q_i}{w_i}}
     $
     \enspace
     and the antecedent of \ref{p7} can be detached
     on account of {\em axiom } (I1).
     We are now supplied with
     \begin{displaymath}
       \dedol{T'}{\smash[b]%%%%%%%%%%%%%%%%%%%%%%%
       \op(\Delta,\Condx va,\Gamma) =
       \op(\Delta,\Condx qa\psubv{v_i}{q_i},\Gamma)
       } %%%%%%%%%%%%%%%%%%%%%%%%%%%%%%%%%%%
     \end{displaymath}
     \noindent and due to \ref{norm.p}(2):
     $\ded{T'}\obl{e_1=e_2} \leqv \eqobl{\nor{e_1}}{\nor{e_2}}$
     we can interchange the two terms (left and right side of equation)
     within the considered context and the requirement for applicability
     of \ref{HS5} will be provided. Hence without loss of generality we can
     assume the required condition for $\vec{x}$ and $\vec{v_i}$, too.
     This confirms (|||) and through the chain
     $(|||) \limp (||) \kom (||) \limp (|)$ we succeed to (|). \qed(3.4)
    To check (4) of \ref{Strukt} Def. that $\CM_{\prop}$
    is a Boolen algebra with two elements, whose operations are the
    $\CM$-values of certain logical connectives: this can be seen from
    \ref{TStr} Def. of {\em term structure } together with
    \ref{norm} Def. of {\em norm} and simple instances of
    propositional tautologies as e.g.
    $\obl{\falsum \lor \verum \eqs{\prop} \verum }$
    (note $\obl{\eqs{\prop}}$ is $\obl{\leftrightarrow}$).
    Finally check (5): From \ref{TStr} Def. $\CM$,
    special case $\oblt{op}=\obl{\eqs{\gamma}}$,
    only to consider case $r_i=0$, we obtain
    $\CM\obl{\eqs{\gamma}}(\obl{\nor{a_1}},\obl{\nor{a_2}}) =
       \obl{a_1 \eqs{\gamma} a_2}
    $
    and have to show
    \begin{math}
       \obl{\nor{a_1 \eqs{\gamma} a_2}} = \obl{\curlyvee}
       \text{ if } \obl{\nor{a_1}} = \obl{\nor{a_2}}
       \text{ otherwise } = \obl{\curlywedge}.
    \end{math}
    This follows from \ref{norm.p}(2)
    $\obl{\nor{a_1}} = \obl{\nor{a_2}}$ iff
    $\T'\dedq \obl{a_1 \eqs{\gamma} a_2}$
    and \ref{norm.p}(1) ... iff
    $\obl{\nor{a_1 \eqs{\gamma} a_2}} = \obl{\curlyvee}$.
 \end{proof}

 \begin{theorem}\em\label{CM.Expr}
    Assume $\CM$ is the {\em term structure } of $\T'$
    due to Def. \ref{TStr},
    $\vec{\sigma} = \enangle{\sigma_i}_{\indto i\ell} \in \VSRT^{\ell}
    \komma
    \vec{x} = \enangle{\obl{x_i}}_{\indto i\ell} \in \VAR_{\vec{\sigma}}
    $, then
 \end{theorem}
 \skipbasel{-1.8}
 \begin{displaymath}
    \sparenth{normal}{\obl{e}\in\Expr_{\gamma}[\vec{x}]} \quad
    \obl{e}\inCM(\vec{x}) =
    \begin{cases}
       \obl{e}\inCM() = \obl{\nor{e}} &\text{ if } \ell = 0
       \\
       \funkd{\CM_{\vec{\sigma}}}{\CM_{\gamma}}%
       {\vec{s}}{\obl{\nor{\sub{e}{\vec{x}}{\vec{s}}}}}
       &\text{ if } \ell \neq 0
    \end{cases}
 \end{displaymath}

 % PROOF72.TEX
 % ~~~~~~~~~~~
 \begin{proof}
    We assume $\obl{e}\in\Expr_{\gamma}[\vec{x}]$ and $\pSyntass$
    (ref. \fullref{pGP}, \ref{TStr} and \fullref{Interp}),
    w.r.t. \ref{Interp} substitute $\obl{e}_{\CM}(\vec{x})$ for
    $\obl{e}_{\M}(\vec{u})$).
    \\
    {\em Case } $\ell=0$: Can be treated like the other case.
    \qquad
    {\em Case } $\ell\neq 0$:
    The goal is to show
    \begin{equation}\tag{|}
       \obl{e}_{\CM}(\vec{x})(\vec{s}) =
       \obl{\nor{e\psubv xs}}
       \quad \text{for } \vec{s} \in \CM_{\sigma}
       \qquad (\text{note that } \CM_{\sigma} \subseteq \Expr_{\sigma}[])
    \end{equation}
    We shall evaluate the left side of the equation according to
    \ref{Interp} Def. of {\em interpretation } and the right side
    due to \ref{subst} and utilize from the {\em induction premise } that the
    law to be shown already applies to the argument terms $\obl{a_i}$
    until {\em left = right } becomes evident.
    {\em Case } $m=0$:
    From $\obl{e}\in\Expr_{\gamma}[\vec{x}]$,
    by \fullref{pGP}, \fullref{vcmps} and
    the supposed $\pSyntass$ we obtain
    $\obl{e}=\oblt{op} \in \SOP \cup \vcmps x
    $.
    First consider case
    $\obl{e}=\oblt{op} \in \vcmps x
    $.
    \eleft
    According to \ref{Interp} (Def. {\em interpretation})
    \begin{math}
       \obl{e}_{\CM}(\vec{x})(\vec{s}) =
       \text{pj}_k^{\vec{\sigma}} = \obl{s_k}.
    \end{math}
    \eright
    \ref{subst} yields $\obl{e\psubv xs}=\obl{s_k}$.
    For both {\em left } and {\em right } $k$ is defined
    to be the biggest so that $\obl{x_k}=\obl{e}$.
    $\obl{s_k} \in \CM_{\sigma_k}$ implies $\obl{s_k}=\obl{\nor{s_k}}$,
    hence $\obl{\nor{e\psubv xs}}=\obl{s_k}$. (|) is evident.
    Now consider the other case $\obl{e}=\oblt{op} \in \SOP$.
    \eleft
    by \ref{Interp}\
    $\obl{e}_{\CM}(\vec{x}) = \text{cst}^{\vec{\sigma}}_{\CM\oblt{op}}$,
    then
    $\obl{e}_{\CM}(\vec{x})(\vec{s}) = \CM\oblt{op}
    \underset{\text{\ref{TStr}(3)}}= \obl{\nor{\op}} = \obl{\nor{e}}$
    \eright by \ref{subst}
    $\obl{\nor{e\psubv xs}} = \obl{\nor{\op}} = \obl{\nor{e}}$

    {\em Case } $m \neq 0$:
    We rely on our {\em induction hypothesis: }
    $\obl{a_i}_{\CM}(\vec{x}\vec{v_i})(\vec{s}\vec{b_i}) =
       \obl{\nor{a_i\psubvv{x}{v_i}{s}{b_i}}}
    $
    and observe $\obl{e} = \obl{\opex va}$.
    \eright $\obl{\nor{e\psubv xs}}$.
    \eleft
    According to \ref{Interp} (Def. {\em interpretation})
    $\obl{e}_{\CM}(\vec{x}) =
     \CM{\oblt{op}}(\enangle{\reob{h}_i \vec{s}}_{\indto im})
    $ with
    $\reob{h}_i\colon\CM_{\vec{\sigma}}\to\CM_{\alpha_i}
    $,
    so that for arbitrary $\vec{s} \in \CM_{\vec{\sigma}}$
    \begin{equation}\notag
       \text{if } r_i=0 \text{ then }
       \reob{h}_i = \obl{a_i}_{\CM}(\vec{x})(\vec{s})
       \qquad \text{else }\quad
       \reob{h}_i =
       \smash{%
        \funkd{\CM_{\vec{\beta_i}}}{\CM_{\alpha_i}}%
       {\vec{b}}{\obl{a_i}_{\CM}(\vec{x}\vec{v_i})(\vec{s}\vec{b})}
       }
    \end{equation}
    \noindent {\em induction hypothesis } yields
    \begin{equation}\tag{i}
       \text{if } r_i=0 \text{ then }
       \reob{h}_i = \obl{\nor{a_i\psubv xs}}
       \qquad \text{else }\quad
       \reob{h}_i =
       \smash{%
        \funkd{\CM_{\vec{\beta_i}}}{\CM_{\alpha_i}}%
       {\vec{b}}{\obl{\nor{a_i\psubvv{x}{v_i}{s}{b}}}}
       }
    \end{equation}
    In order to apply \ref{TStr}(Def. $\CM$)
    we transform the equation next to label ``\eleft'' into
    \begin{math}
       \opaar%
       {\enangle{\reob h_i \vec{s}}_{\indto im}}%
       {\obl{e}_{\CM}(\vec{x})(\vec{s})}
       \in \CM\oblt{op}
    \end{math}
    and in order to distinguish already used symbols (meta variables)
    from corresponding ones that may denote different objects
    we supply the one with star as superscript.
    Applying \ref{TStr}(3) in this way we obtain
    $\CM\oblt{op}\colon\Prodim\CM_{\alpha_i}^{\vec{\beta_i}}
       \to \CM_{\gamma}
    $
    and that there are
    $\vec{\avar{a}} \kom \VVec{\avar{v}} \kom \obl{\avar{e}}$
    so that
    \begin{math}
       \VVec{\avar{v}} = \enangle{\vec{\avar{v_i}}}_{\indto im}
       \komma
       \vec{\avar{a}} =
       \enangle{\obl{\avar{a_i}}}_{\indto im}
       \komma
       \text{each } \vec{\avar{v_i}} \in \VAR_{\vec{\beta_i}}
       \komma
       \text{each }
       \obl{\avar{a_i}} \in \Expr_{\alpha_i}[\vec{\avar{v_i}}]
    \end{math}
    and (as $\enangle{\reob{h}_i}_{i}$ substitutes $\vec{\frak{a}}$)
    \vspace{-1ex}
    \begin{equation}\tag{ii}
       \text{if } r_i=0 \text{ then }
          \reob{h}_i(\vec{s}) = \obl{\nor{\avar{a_i}}}
       \text{ else }
          \reob{h}_i(\vec{s}) =
          \funkd{\CM_{\vec{\beta_i}}}{\CM_{\alpha_i}}{\vec{b}}%
             {\obl{\nor{\avar{a_i}\psubv{\avar{v_i}}{b}}}}
       \komma
    \end{equation}
    finally
    \begin{math}
       \obl{\nor{\avar{e}}} = \obl{\nor{\opex{\avar{v}}{\avar{a}}}}
    \end{math}
    and
    \begin{math}
       \obl{e}_{\CM}(\vec{x})(\vec{s}) = \obl{\nor{\avar{e}}}.
    \end{math}
    Note that each starred metavariable depends on parameter $\vec{s}$.
    Our {\em goal } is to deduce
    $\obl{\nor{\avar{e}}}=\obl{\nor{e\psubv xs}}$.
    To achieve it we try to obtain the premises suited for an
    application of \ref{eqlemma2} {\em lemma}.
    Combining (i)+(ii) yields
    \begin{equation}\notag %%% {iii}
       \text{if } r_i=0 \text{ then }
       \obl{\nor{a_i\psubv xs}} = \obl{\nor{\avar{a}_i}}
       \text{ else }
       \AQu{\vec{b}\in\CM_{\vec{\beta_i}}}
       \obl{\nor{a_i\psubvv{x}{v_i}{s}{b}}}
       = \obl{\nor{\avar{a_i}\psubv{\avar{v_i}}{b}}}
    \end{equation}
    Note that case $r_i=0$ can be treated like the other case,
    as $\obl{\avar{a}_i[\emptyseq\gets\emptyseq]}=\obl{\avar{a}_i}$
    is stipulated by \ref{subst} {\em def.},
    \ref{norm.p} yields
    \begin{math}
       \T'\dedq \obl{%
          a_i\psubvv{x}{v_i}{s}{b} =
          \avar{a}_i\psubv{\avar{v_i}}{b}
       }
    \end{math},
    \\
    by \ref{HS3} \enspace
    \begin{math}
       \obl{a_i\psubvv{x}{v_i}{s}{b}} =
       \obl{a_i\psubvv{x}{v_i}{s}{v_i}\psubv{v_i}{b}}
    \end{math}
    \enspace hence
    \begin{equation}\tag{iii}
       \AQu{\vec{b}\in\CM_{\beta}} \quad
       \T'\dedqol{a_i\psubvv{x}{v_i}{s}{v_i}\psubv{v_i}{b}
          = \avar{a}_i\psubv{\avar{v_i}}{b}
          }
    \end{equation}
    This is the prerequisite for applying \ref{eqlemma2} which yields
    the required $\obl{\nor{\avar{e}}}=\obl{\nor{e\psubv xs}}$.
 \end{proof}

 \begin{observ}\em\label{CM.verum}
    \begin{math}
       \CM\obl{\verum} = \obl{\verum} = \AlgOp1
       \qquad
       \CM\obl{\falsum} = \obl{\falsum} = \AlgOp0
    \end{math}
 \end{observ}

 \begin{observ}\em
    If $\T\subseteq\Expr_{\prop} \komma \obl{\varphi} \in \Expr_{\prop}$
    then $\T\dedqol{\AQ{\vec{x}}\,\varphi} \leqv \T\dedqol{\varphi}$
 \end{observ}

 \begin{theorem}\em\label{ded=sat}
    Assume the premises as in \ref{CM.Expr}, let
    $\gamma=\prop$ and $\obl{\varphi} \in \Expr_{\prop}[\vec{x}]$, then
 \end{theorem}
 \EquaNull
 \vspace{-0.7\baselineskip}
 \begin{alignat}{2}
    \T'\dedq \obl{\varphi}  &\leqv \obl{\nor{\varphi} }\inCM()
     = \obl{\verum} = \AlgOp1
    &&\qquad\text{if } n = 0
    \\
    \T'\dedq \obl{\varphi}
    &\leqv
    \AQu{\vec{a} \in \CM_{\vec{\sigma}}}\enspace
       \obl{\nor{\varphi}}\inCM(\vec{x})(\vec{a})
        = \obl{\verum} = \AlgOp1
    &&\qquad\text{if } n \neq 0
    \\
    \T'\dedq \obl{\varphi}  &\leqv
    \CM \sat{\rule{0pt}{1.8ex} S'} \obl{\varphi}
    &&\qquad(n \ge 0)
 \end{alignat}
 \begin{proof}
    The two preceding \textit{observations} and application of
    \ref{CM.Expr}  and \ref{norm.p} yield this result.
 \end{proof}

 \begin{observ}\noem\label{sat.restr}
    Let $\T$ be a consistent theory and
    $\HT$ its extension due to \fullref{ErwThm}
    ({\em extension theorem}).
    If $\CM$ is a model of $\T'$ then the restriction
    $S\restr\CM$ of $\CM$ to the {\em signature }
    $S \sqsubseteq S'$ is a model of $\T$, that is:
    \quad
    $\CM \sat{S'} \T' \limp S\restr\CM \sat{S} \T$
 \end{observ}

 \begin{thms}{Completeness Theorem (2nd Version)}
    \noem\label{vollst2}
    A consistent {$\fnl-$\em theory } has a model.
 \end{thms}
 \skipbasel{-1}
 \begin{displaymath}
    \hspace*{-10mm}
    \fnlsig S \land \T \subseteq \Expr_{\prop}
    \land  \kons(\T)
    \limp
    \EQu{\M \, \FnlSMG}\,  \M \sat{S} \T
 \end{displaymath}

 \begin{pf}
    This is a consequence of
    \ref{ErwThm},\ref{TStr},\ref{ded=sat},\ref{sat.restr}
    according to \ref{sat}.
 \end{pf}

 \begin{thms}{Completeness Theorem (G\"odel)}\noem
    Any formula that is valid in a {$\fnl-$\em theory }
    is provable in it, that is
 \begin{math}
    \fnlsig S \land
    \T \cup \encurs{\obl{\varphi}} \subseteq \Expr_{\prop}
    \land \T \sat{S} \obl{\varphi} \limp \T \ded{S} \obl{\varphi}
 \end{math}
 \end{thms}

 \begin{proof}^{\newcommand{\TN}{\T \cup \encurs{\obl{\neg\varphi}}}}
    \begin{tabular}[t]{R@{\quad iff\quad}L}
       \T \sat{S} \obl{\varphi}
       &\TN
       \hspace{1.5em} \text{has no model}
       \\
       \T \ded{S} \obl{\varphi}
       &\mathbf{\lnot} \kons(\TN)
    \end{tabular}
    \\[3pt]
    Use \ref{vollst2} with $\TN$ in place of $\T$,
    then the first 2 premises of the theorem imply
    \\[3pt]
    \begin{mArray}^{\renewcommand{\arraystretch}{1.3}}{r@{\limp}l}
       \kons(\TN) &\EQu{\M}\,\M \sat{S} \TN
       \\
       \mathbf{\lnot} \EQu{\M}\,\M \sat{S} \TN
       &\mathbf{\lnot} \kons(\TN)
       \qquad
       \text{by \it propositional tautology}
       \\
       \T \sat{S} \obl{\varphi}
       &\T \ded{S} \obl{\varphi}
       \qquad
       \text{according to the first two lines of the proof}
    \end{mArray}
 \end{proof}

 \noindent\textbf{Remarks.}
 The construction of a \textit{term model} could be carried out also
 with a theory $\T'$, if the requirement of \vref{Henkin}
 (or \vref{ErwThm} respectively) is weakened
 by rewriting $\COP_{S'}^{\gamma}$ into $\Expr_{S'}^{\gamma}[]$
 (note that $\COP_{S'}^{\gamma} \subseteq \Expr_{S'}^{\gamma}[]$).
 We then say $\T'$ \textit{admits examples}.
 If we prove the validity of
 \textit{Hilbert's 2nd $\varepsilon$-Theorem} % (ref. to \cite{epsilon})
 for \textit{Functional Logic}, which claims
 $\T'$ to be a \textit{conservative extension} of $\T$,
 if the new symbols are $\varepsilon_{\gamma}$ with signature
 \begin{math}
    \obl{\gamma((\gamma):\prop)}
 \end{math}
 for each $\gamma\in\VSRT$ and the additional \textit{axioms}
 are $(\varepsilon_0)$ and $(\varepsilon_2)$ as defined in \cite{epsilon}
 (but actually multiplied by indexing range $\VSRT$),
 then $\T'$ \textit{admits examples} and
 inherits \textit{consistency} from $\T$.
 We observe that
 \begin{monolist}[\labelwidth=1.6em\labelsep=.7em]{($\varepsilon_2$)}
    \begin{math}
       \AQ{z}(\varphi^x_z \leqv \psi^y_z)
       \limp \varepsilon x.\varphi = \varepsilon y.\psi
    \end{math}
    \quad
    is accurately the \textit{equality ax.} (I2)
    for the symbol $\obl{\varepsilon}=\obl{\varepsilon_{\xi}}$
    (see p. \pageref{CAL} and remember that
    $\obl{\leftrightarrow}$ is the same as $\obl{\eqs{\prop}}$).
   \item[($\varepsilon_0$)]
       $\EQ{x}\,\varphi \limp \varphi\psub{x}{\varepsilon x.\varphi}$
       \quad
       is a nonlogical axiom, which makes $\T'$ admit examples.
 \end{monolist}
 Thus $\T'$ with only one new symbol for each $\gamma \in \VSRT$
 would be equally suitable for constructing a term model. Yet proving
 \textit{Hilbert's 2nd $\varepsilon$-Theorem} (ref. to \cite{epsilon})
 for this purpose
 is considerably more difficult and largescale than Henkins approach.
 \endgroup %%%<<<<<<<<<<<<<<<<

 \vfill \small
%
%%####################################################################
%%# funlk.bbl   |95-03-12|09:50|                                     #
%%####################################################################
 